\documentclass[12pt]{article}
\usepackage{amsmath,amssymb}
\def\MC{\mathcal{M}}

\def\LC{\mathcal{L}}
\def\XC{\mathcal{X}}
\def\CC{\mathcal{C}}
\def\UC{\mathcal{U}}
\def\FC{\mathcal{F}}

\def\R{\mathbf{R}}
\def\N{\mathbf{N}}
\def\x{\mathbf{x}}
\def\y{\mathbf{y}}
\def\E{\mathbf{E}}
\def\1{\mathbf{1}}

\newtheorem{proposition}{Proposition}[section]
\newtheorem{theorem}{Theorem}[section]
\newtheorem{lemma}{Lemma}[section]

\numberwithin{equation}{section}

\begin{document}
\title{THE CENTRAL LIMIT THEOREM FOR THE SMOLUCHOVSKI COAGULATION MODEL}
\author{Vassili N. Kolokoltsov\thanks{Department of Statistics, University of Warwick, Coventry CV4 7AL, UK.
 Email: v.kolokoltsov@warwick.ac.uk}}
\date{}
\maketitle
\begin{abstract} \it
The general model of coagulation is considered.
 For basic classes of unbounded coagulation kernels the central limit
  theorem (CLT) is obtained for the fluctuations around the dynamic law
   of large numbers (LLN). A rather precise rate of convergence is given both for LLN and CLT.
\end{abstract}

\section{Introduction}

Throughout the paper we shall denote by $X$ a locally compact topological
 space equipped with its Borel sigma algebra and by $E$ a given continuous
  non-negative function on $X$ such that $E(x) \rightarrow \infty$ as $x \rightarrow \infty$ .
   Denoting by $X^0$ a one-point space and by $X^j$ the powers
 $X \times \ldots \times X $($j$-times) considered with their product topologies,
  we shall denote by $\XC$ their disjoint union $\XC = \bigcup^\infty_{j=0}X^j$,
   which is again a locally compact space. In applications, $X$ specifies
  the state space of a single particle, $\XC$ stands for the state space
  of a random number of similar particles, and $E$ describes some key parameter
   of a particle. In the standard model $X = \R_+ = \{x > 0\}$ and $E(x) = x$ denotes the mass of a particle.

By $C(X)$ (respectively $C_\infty(X)$) we always denote the Banach
space of continuous bounded functions on $X$ (respectively its
subspace of functions vanishing at infinity) with the sup-norm
denoted by $\|\cdot\|$, by $\MC(X)$ - the Banach space of finite
Borel measures on $X$  with the norm also denoted by $\|\cdot\|$,
and by $\MC^+(X)$ - the set of its positive elements. The brackets
$(f,Y)$ denote the usual pairing (given by the integration) between
functions $f$ and measures $Y$, and $|\mu|$ for a signed measure
$\mu$ denotes its total variation measure. The elements of $\XC$
will be denoted by bold letters, e.g. $\x = (x_1, \ldots, x_n) \in
X^n \subset \XC$. For a subset $I$ in $\{1, \ldots, n\}$ we shall
denote by $|I|$ and $\bar{I}$ respectively its cardinality and its
complement
 in $\{1, \ldots, n\}$, and by $\x_I$ the element of $X^{|I|}$ given by the collection of $x_i$ with $i \in I$.

Assume that we are given a continuous transition kernel $K(x_1,
x_2;dy)$ from $X \times X$ to $X$, i.e. a continuous function from
$X \times X$ to $\MC^+(X)$ (the latter equipped with its $*$-weak
topology, i.e. the topology of the dual space to $C_{\infty}(X)$).
This kernel will be called the coagulation kernel and it will be
assumed to preserve $E$, i.e. $K(x_1, x_2;dy)$ has support contained
in the set $\{y: E(y) = E(x_1) + E(x_2)\}$. Moreover, $K(x_1,
x_2;dy)$ is symmetric with respect to permutation of $x_1$ and $x_2$
and has intensity $K(x_1,x_2) = \int_X K(x_1,x_2;dy)$ enjoying the
following additive upper bound:
\begin{equation}
\label{eq1.1} K(x_1,x_2) \leq C(1 + E(x_1) + E(x_2))
\end{equation}
with some constant $C > 0$ and all $x_1,x_2$.

The process of coagulation that we are going to analyse here is
a Markov process $Z(t)$ on $\XC$ specified by the generator
\begin{equation}
\label{eq1.2} Lg(\x) = \sum_{I \subset \{1, \ldots, n\}:|I|=2} \int
(g(\x_{\bar{I}},y) - g(\x))K(\x_I ; dy)
\end{equation}
(where $\x = (x_1, \ldots, x_n)$) of its Markov semigroup acting on
an appropriate space
 of functions on $\XC$. It is known and not difficult to deduce from the theory of jump
 type processes (see e.g. \cite{Ch}) that the process $Z(t)$ is well defined by this generator
  (see e.g. a detailed probabilistic description of $Z(t)$ in \cite{No}). In the next
  Section the analytic properties of the Markov semigroup specified by $L$ will be made precise.

The transformation

\begin{equation}
\label{eq1.3}
\x = (x_1, \ldots, x_n) \mapsto h \delta_{\x} = h(\delta_{x_1} + \cdots \delta_{x_n}),
\end{equation}
with $h$ being a positive (scaling) parameter, maps $\XC$ to the space $\MC_{h \delta}(X)$
of positive measures on $X$ of the form $h\delta_{\x}$. By $Z^h_t$ we shall denote
a Markov process on $\MC_{h \delta}(X)$ obtained from $Z(t)$ by transformation
\eqref{eq1.3} combined with the scaling of $L$ by $h$, i.e. $Z^h_t$ is defined through the generator
\begin{align}
\label{eq1.4} & L_h G_g(h \delta_{\x})
= h \sum_{I \subset \{1, \ldots, n\}:|I| = 2} \int (g(\x_{\bar{I}},y) - g(\x))K(\x_I;dy) \nonumber \\
& = h \sum_{I \subset \{1, \ldots, n\}:|I| = 2} \int (G_g(h\delta_\x + h(\delta_y - \delta_{\x_I}))- G_g(h \delta_\x)) K(\x_I;dy)
\end{align}
on $C(\MC_{h \delta}(X))$, where $G_g(h \delta_\y) = g(\y)$ for any
$\y \in \XC$.

The law of large numbers dynamics (LLN) for the processes $Z^h_t$ is
given by the kinetic equation,
 whose most natural form is the weak one, i.e. it is the equation
\begin{equation}
\label{eq1.5}
\frac{d}{dt} (g, \mu_t)
= \frac{1}{2} \int_{X \times X} \int_X (g(y) - g(x_1) - g(x_2)) K(x_1,x_2;dy) \mu_t(dx_1) \mu_t (dx_2)
\end{equation}
on $\mu_t$ that has to hold for all $g \in C_\infty(X)$. It is known (see \cite{No})
 that if a family of initial measures $h \delta_{\x(h)}$ for $Z^h_t$ is uniformly
  bounded with bounded moments of order $\beta \geq 2$, i.e. if
\begin{equation}\label{eq1.6}
\sup_h \int_X (1 + E^\beta(y)) h \delta_{\x(h)}(dy) < \infty,
\end{equation}
and if $h \delta_{\x(h)}$ tends $*$-weakly to a measure $\mu_0$ on $X$, as $h \rightarrow 0$,
 then the process $Z^h_t$ with the initial data $h \delta_{\x(h)}$ tends weakly to
 a bounded solution $\mu_t$ of \eqref{eq1.5} with initial condition $\mu_0$ that preserves
  $E$ and has bounded moments of order $\beta$, i.e. such that
\begin{equation}\label{eq1.7}
\sup_{s \leq t} \int_X (1 + E^\beta(y))\mu_s(dy) < \infty
\end{equation}
and
\begin{equation}\label{eq1.8}
\int_X E(y) \mu_t (dy) = \int_X E(y) \mu_0(dy)
\end{equation}
for all $t \geq 0$.

The first objective of this paper is to establish the corresponding central limit theorem (CLT),
 i.e. to show that the process
$$F^h_t(Z^h_0, \mu_0) = h^{-1/2} (Z^h_t (Z^h_0) - \mu_t(\mu_0))$$
of normalized fluctuations of $Z^h_t$ around its dynamic law of
large numbers
 $\mu_t$ converges in some sense to a generalized Gaussian Ornstein-Uhlenbeck process on
 $\MC(X)$ or a more general space of distributions. We obtain this result under some mild
  technical assumptions on the coagulation kernel thus presenting a solution to the problem 10
   from the list of open problems on coagulation formulated in the well known review \cite{Al}.

It is worth noting that though for the classical processes preserving the number of particles
 (like interacting diffusions or Boltzmann type collisions) the results of CLT type are well
  established and widely presented in the literature (see e.g. \cite{GW} or \cite{Da1}
   and references therein), for the processes with a random number of particles the work on CLT began recently.
For coagulation processes with discrete state space $X = \N$ and
uniformly bounded intensities the central limit for fluctuations was
obtained in \cite{DFT} using stochastic
 calculus. For general processes of coagulation, fragmentation and collisions on $X = \R_+$,
  but again with bounded intensities, the central limit was proved by a different method in \cite{Ko6},
   namely by analytic methods of the theory of semigroups. The results of the present paper are obtained
    by developing further the approach from \cite{Ko6}.

The second objective of the paper is to provide precise estimates of the error term both
 in LLN and CLT for a wide class of bounded and unbounded functionals on measures. Note that
  the usual ``prove compactness in the Skorohod space and choose a converging subsequence''
  probabilistic method does not provide such estimates (see, however, \cite{GW} for a progress
  in this direction for interacting diffusions). Our main technical tool is the study of the
   derivatives of the solutions to kinetic equations with respect to initial data (this approach
    is inspired by the analysis of such derivatives for the Boltzmann equation in \cite{Ko4}).
     The existence and regularity of these derivatives in weighted spaces of functions and measures
      are analyzed and the validity of CLT is proved to be connected with a certain kind of
      stability of these derivatives. The final estimates and their proofs depend on the structure
   and the regularity properties of the coagulation kernel. We demonstrate various aspects of our
    approach analyzing the following three classes of kernels:

 (C1) $K(x_1,x_2) = C(E(x_1) + E(x_2))$.

{\em Remark.} This is a warming up example, for the solutions to the main equations
 are given more or less explicitly in this case.

 (C2) $K(x_1,x_2) \leq C(1 + \sqrt{E(x_1)})(1 + \sqrt{E(x_2)})$.

 {\em Remark.} This model is analyzed to show the kind of results one can expect to
 obtain without assuming any differential or linear structure on the state space $X$.
  The unavoidable shortcoming of these results is connected with the absence of an appropriate
   space of generalized functions to work with. Hence the estimate of errors in LLN and CLT
   have to depend on something like the norm of $F^h_0 = (Z^h_0 - \mu_0) / \sqrt{h}$ in $\MC(X)$.
    But general $\mu_0$ can not be approximated by Dirac measures $Z^h_0$ in such a way that
     $F^h_0$ be bounded in $\MC(X)$. Hence the possibility to apply these results beyond discrete
      supported initial measures $\mu_0$ is rather reduced. On the other hand, these kind of results
      are open to extensions to very general spaces.

 (C3) $ X = \R_+$,
  $K(x_1, x_2, dy) = K(x_1, x_2) \delta (y - x_1 - x_2)$, $E(x) = x$, $K$ is
  non-decreasing in each argument 2-times
  continuously differentiable on $(\R_+)^2$ up to the boundary
  with all the first and second partial derivatives being bounded by a constant $C$.

  {\em Remark.}  This is the case of our main interest. Unlike previous cases the estimate here
  turns out to depend on the norm of $F^h_0$ coming from the dual space to continuously
  differentiable functions, and this norm can be easily made small for an arbitrary measure
   $\mu_0$ on $X$. Therefore, to shorten the exposition, we shall
   prove  CLT completely, up to the convergence of the distributions
   of processes on the Skorohod space of c\`adl\`ag
   functions, only for this case, restricting the discussion of the first two cases only
   to the convergence of linear functionals.
   For simplicity, we choose here the state space $X=\R_+$
    of the standard Smoluchowski model, the extensions to finite-dimensional Euclidean spaces $X$ being not difficult
    to obtain. Similarly we choose very strong assumptions on the derivatives (in particular,
     the kernels $K(x,y) = x^\alpha + y^\alpha$ with $\alpha \in (0,1)$
    are excluded by our assumption, as the derivatives of this $K$ have a singularity
    at the origin). Finally let us stress that all kernels from (C1)-(C3) clearly satisfy
     \eqref{eq1.1} (possibly up to a constant multiplier).

  We refer to reviews \cite{Al} and \cite{Le} for a general background in coagulation models,
   and to \cite{KSW} for simulation and numerical methods.

  The content of the paper is the following. In the next section we formulate the main results,
   and other sections are devoted to their proofs. In particular, Sections 4 and 5 are devoted
    to a detailed analysis of the equation in variations (linear approximation)
     around the solution of kinetic equation \eqref{eq1.5} that describes the derivatives
      of the solution to \eqref{eq1.5} with respect to the initial measure $\mu_0$.
       At the end of Sect. 5 a new property of the kinetic equation itself is
       established that is crucial to our proof of CLT, but seems to be also of
       independent interest.
  Namely Propositions \ref{prop5.5}, \ref{prop5.8} show that the solution depends Lipschitz continuously on the
  initial measure in the topology of the dual to the weighted spaces of continuously differentiable
   functions or certain weighted Sobolev spaces. In Section 9 three
   general result are presented (on variational derivatives, on the
   linear transformation of Feller processes and on the dynamics of
   total variations of measures), used in our proofs and places separately in order
   not to interrupt the main line of arguments.
In the Appendix some auxiliary facts  on the evolutions specified by
unbounded integral generators are presented. Though they should be
essentially known to probabilists dealing with jump processes, the
author did not find an appropriate reference.

  To conclude the introduction we shall fix the basic notations to be used throughout the paper
   without further reminder recalling as we go some relevant facts
   about Sobolev spaces and variational derivatives.

  (i) {\it Weighted spaces of functions and measures arising from unbounded intensities of jumps.}
   For a positive measurable function $f$ on a topological space $T$ we denote by
    $C_f = C_f(T)$ and $B_f = B_f(T)$ (omitting $T$ when no ambiguity may arise) the Banach
    spaces of continuous and measurable functions on $T$ respectively having finite norm
  $$
  \| \phi\|_f = \| \phi\|_{C_f(T)} = \sup_x (| \phi(x)| / f(x)).
  $$
  By $C_{f, \infty} = C_{f, \infty}(T)$ and $B_{f, \infty} = B_{f, \infty}(T)$ we denote
  the subspaces of $C_f$ and $B_f$ respectively consisting of functions $\phi$ such that
   $(\phi / f)(x) \rightarrow 0$ as $f(x) \rightarrow \infty$. If $f$ is a continuous function
    on a locally compact space $X$ such that $f(x) \rightarrow \infty$, as $x \rightarrow \infty$,
     then the dual space to $C_{f, \infty}(X)$ is given by the space $\MC_f(X)$ of Radon measures
      on $X$ with the norm $\|Y\|_f = \sup \{ (\phi, Y) : \| \phi\|_f \leq 1\}$.

We shall need also the weighted $L_p$ spaces. Namely, define
$L_{p,f}=L_{p,f}(T)$ as the space of measurable functions $g$ on a
measurable space $T$ having finite norm $\Vert
g\Vert_{L_{p,f}}=\Vert g/f\Vert_{L_p}$.

  For $X= \R_+=\{x>0\}$ we shall use also smooth functions. For a positive $f$ we denote by $C^{1,0}_f(X)$
   the Banach space of continuously differentiable functions $\phi$ on $X = \R_+$ such that
   $\lim_{x \rightarrow 0} \phi(x) = 0$ and the norm
$$
\| \phi\|_{C^{1,0}_f(X)} = \| \phi'\|_{C_f(X)}
$$
is finite. By $C^{2,0}_f(X)$ we denote the space of two-times
continuously differentiable
 functions such that $\lim_{x \rightarrow 0} \phi(x) = 0$ and the norm
$$
\| \phi\|_{C^{2,0}_f(X)} = \| \phi'\|_f  + \| \phi''\|_f
$$
is finite. By $\MC^1_f(X)$ and $ \MC^2_f(X)$ we shall denote the
Banach dual spaces to $C^{1,0}_f$ and $C^{2,0}_f$
 respectively. Actually we need only the topology they induce on (signed) measures so that for
  $\nu \in \MC(X)\cap \MC^i_f(X)$, $i=1,2$,
\[
\| \nu \|_{\MC^i_f(X)}  = \sup \{ (\phi, \nu) : \|
\phi\|_{C^{i,0}_f(X)} \leq 1 \}.
\]

Similarly one defines the spaces $L_{p,f}^{1,0}$ and
$L_{p,f}^{2,0}$, $p\ge 1$, as the spaces of absolutely continuous
functions $\phi$ on $X = \R_+$ such that
   $\lim_{x \rightarrow 0} \phi(x) = 0$ with the norms respectively
$$
\| \phi\|_{L^{1,0}_{p,f}(X)} = \| \phi'\|_{L_{p,f}(X)}
 = \|\phi'/f\|_{L_p(X)},
 \quad
 \| \phi\|_{L^{2,0}_{p,f}(X)} = \| \phi'/f\|_{L_p(X)}
  +\|(\phi'/f)'\|_{L_p(X)},
$$
as well as their dual $(L_{p,f}^{1,0})'$ and $(L_{p,f}^{2,0})'$.

  As an example (needed later) let us
  estimate two of these norms for the Dirac measure $\delta_x$ on $\R_+$, $x>0$ and the function
  $f(y)=f_k(y)=1+y^k$:
  $$
  \|\delta_x\|_{\MC^1_{f_k}(\R_+)}
  =\sup \{ \int_0^x g(y)\, dy: \|g\|_{C_{f_k}}\le 1\}
  =x+x^{k+1}/(k+1);
  $$
 \begin{equation}\label{eq1.81}
 \|\delta_x\|_{(L^{1,0}_{2,f_k})'(\R_+)}
  =\sup \{ \int_0^x g(y)\, dy: \|g/f_k\|_{L_2}\le 1\}
  \le \sqrt{\int_0^x f_k^2 (y)\, dy}\le c(k)\sqrt xf_k(x).
  \end{equation}

 Not every $\nu \in \MC(X)$ belongs to
$\MC^1_f(X)$ or $ \MC^2_f(X)$. Suppose that $f$ is non-decreasing
and $\nu \in \MC(X)$ is such that
\begin{equation}\label{eq1.811}
\tilde \nu (x)=\int_x^{\infty} \nu (dy)=o(1) (xf(x))^{-1}, \quad
x\to \infty.
\end{equation}
  Then by integration by parts for $g\in
C_f^{1,0}(\R_+)$
$$
(g,\nu)=-\int_0^{\infty} g(x) d\tilde \nu (x)
 =\int_0^{\infty} g'(x) \tilde \nu (x) dx
 $$
 (the boundary term vanish by \eqref{eq1.811}),
 so that
 $$
 \| \nu \|_{\MC^1_f(X)}=\| \tilde \nu \|_{L_{1,1/f}}
 $$
 and
 $$
 \| \nu \|_{\MC^2_f(X)}=\sup \{ (\phi,\tilde \nu): \|
 \phi\|_{C_f}+\|\phi'\|_{C_f} \le 1\}.
 $$
Similarly, as
 $$
 \int_0^x \phi (s)ds \le \|\phi\|_{L_{p,f}} \left( \int_0^x
 f^q(y)dy\right)^{1/q}, \quad \frac {1}{p}+\frac{1}{q}=1,
 $$
 it follows that if $\nu \in \MC(X)$ is such that
 $$
 \tilde \nu =o(1) \left(\int_0^xf^q (y)dy\right)^{-1/q}, \quad x\to
 \infty,
 $$
 then
 $$
 \| \nu \|_{(L^{1,0}_{p,f})'}
  =\| \tilde \nu \|_{L_{q,1/f}}, \quad \frac {1}{p}+\frac{1}{q}=1,
 $$
 $$
 \| \nu \|_{(L^{2,0}_{p,f})'}
  =\sup \{ (\psi,\tilde \nu): \|\psi /f\|_{L_p}+\|(\psi /f)'\|_{L_p} \le
  1\}, \quad p>1.
 $$
 In particular, recalling that the usual Sobolev Hilbert spaces
 $H^k(\R)$ are defined as the completion of the Schwarz space
 $S(\R)$ with respect to the scalar product
 $$
 (f,g)_{H^k}=(f, (1-\Delta)^kg)_{L_2}
  =(\FC f, (1+p^2)^k\FC g)_{L_2},
 $$
 where
 $$
 (\FC(f))(p)=(2\pi)^{-1/2}\int_{\R} e^{-ipx}f(x)\, dx
 $$
 denotes the usual Fourier transform, and that by duality
 $(H^k)'=H^{-k}$ it follows that
$$
 \| \nu \|_{(L^{2,0}_{2,f})'}
  =\sup \{ (\psi,\tilde \nu): \|\psi /f\|_{H^1} \le 1\}
  =\sup \{ (\phi,f\tilde \nu): \|\phi \|_{H^1} \le 1\}
 $$
\begin{equation}\label{eq1.82}
  =\| f\tilde \nu \|_{H^{-1}}
  =\sqrt {\int_{-\infty}^{\infty} |\FC(f\tilde \nu)(p)|^2
  \frac{dp}{1+p^2}}.
  \end{equation}
We shall use this formula in Section 7.

 (ii) {\it Functional spaces describing indistinguishable
particles.} By $C^{sym}(X^k)$ we denote the Banach space of
symmetric (with respect to all permutations of its arguments)
continuous bounded functions on $X^k$, and by $C^{sym}(\XC)$- the
Banach space of continuous bounded functions on $\XC$) whose
restrictions on each $X^k$ belong to $C^{sym}(X^k)$. For a function
$f$ on $X$ we denote by $f^\otimes$ its natural lifting on $\XC$,
i.e. $f^\otimes(x_1, \ldots, x_n) = f(x_1) \cdots f(x_n)$.

If $f$ is a positive function on $X^m = \R_+^m$, we denote by $C^{1, sym}_f(X^m)$
 (respectively $C^{2,sym}_f(X^m)$) the space of symmetric continuous differentiable
 functions $g$ on $X^m$ (respectively two-times continuously differentiable) vanishing
  whenever at least one argument vanishes, with the norm
$$
\|g\|_{C^{1,sym}_f(X^m)}
 = \left\| \frac{\partial g}{\partial x_1} \right\|_{C_f(X^m)}
  = \sup_{x,j} \left( \left| \frac{\partial g}{\partial x_j} \right| (f^{-1}) \right) (x)
  $$
and respectively
$$
\|g\|_{C^{2,sym}_f(X^m)}
= \left\| \frac{\partial g}{\partial x_1} \right\|_{C_f(X^m)}
 + \left\| \frac{\partial^2 g}{\partial x^2_1} \right\|_{C_f(X^m)}
  + \left\| \frac{\partial^2 g}{\partial x_1 \partial x_2} \right\|_{C_f(X^m)} .
  $$

(iii) {\it Variational derivatives.} For a function $F$ on $\MC_f(X)$ the variational derivative
$\delta F$ is defined by
$$
\delta F(Y;x) = \lim_{s \rightarrow 0_+} \frac{1}{s} (F(Y + s \delta_x) - F(Y)),
$$
where $\lim_{s \rightarrow 0_+}$ means the limit over positive $s$.
Occasionally we shall omit the last argument here writing $\delta
F(Y)$ instead of $\delta F(Y;.)$. The higher derivatives $\delta^l
F(Y;x_1,...,x_l)$ are defined inductively.

As it follows from the definition, if $\delta F(Y;.)$ exists and
depends continuously on $Y$ in the $\star$-weak topology of $\MC$
(or any $\MC_f$), then the function $F(Y+s\delta_x)$ of $s\in
{\R}_+$ has a continuous right derivative everywhere and hence is
continuously differentiable, which implies that
 \begin{equation}
\label{eq1.83}
 F(Y+\delta_x)-F(Y) =\int_0^1 \delta F(Y+s \delta_x;x)\, ds.
 \end{equation}

 We shall need an extension of this identity for more general measures in the place of the Dirac measure
 $\delta_x$. To this end
the following definitions turn out to be useful.
 For two continuous
functions $\phi, f$ such that $0 \leq \phi \leq f$ and $f(x)
\rightarrow \infty$ as $x\to \infty$, we say that $F$ belongs to
$C^l(\MC_{f,\phi}(X))$, $l = 0,2, \ldots,$ if $F\in C(\MC_f)$ and
for all $k=1,...,l$, $\delta^kF(Y;x_1, \ldots, x_k)$ exists for all
$x_1, \ldots, x_k \in X^k$, $Y \in \MC_f(X)$ and represents a
continuous mapping $\MC_f(X) \mapsto C_{\phi \otimes \cdots \otimes
\phi,\infty}^{sym}(X^k)$, where $\MC_f(X)$ is considered in its
$*$-weak topology. We shall write shortly $C^l(\MC_f(X))$ for
$C^l(\MC_{f,f}(X))$. All necessary formulae on the variational
derivatives in these classes are collected in Lemma \ref{lem9.1}.

 {\it Remark.} The introduction of the cumbersome notations
$C^m(\MC_{f, \phi}(X))$ is motivated by the fact that (under our
assumption on the growth of the coagulation rates) if one considers
the solution to the kinetic equations $\mu_t$ with $\mu_0 \in
\MC_{1+E^{\beta}}$, then usually $\dot \mu_t \in
\MC_{1+E^{\beta-1}}$ and the derivatives of $\mu_t$ with respect to
the initial data belong to $C_{1+E^k}$ with certain $k<\beta$, see
Sections 4 and 5.

(iv) {\it Propagators.} If $S_t$ is a family of topological linear spaces,
$t \in \R^+$, we shall say that a family of continuous linear operators
$U^{t,r} : S^r \mapsto S^t$, $r \leq t$ (respectively $t \leq r$) is a propagator
(respectively a backward propagator), if $U^{t,t}$ is the identity operator in $S^t$
for all $t$ and the following propagator equation (called Chapman-Kolmogorov equation
 in the probabilistic context) holds for $r \leq s \leq t$ (respectively for $t \leq s \leq r$):
\begin{equation}\label{eq1.12}
U^{t,s} U^{s,r} = U^{t,r}.
\end{equation}

By $c$ and $\kappa$ we shall denote various constants indicating in brackets
 (when appropriate) the parameters on which they depend.

For an operator $U$ in a Banach space $B$ we shall denote
by $\|U\|_B$ the norm of $U$ as a bounded linear operator in $B$.

At last, we shall use occasionally the obvious formula
\begin{equation}
\label{eq1.13} \sum_{I \subset \{1, \ldots, n\}, |I| = 2} f(\x_I) =
\frac{1}{2} \int \int f(z_1, z_2) \delta_\x(dz_1)\delta_\x(dz_2) -
\frac{1}{2} \int f(z,z) \delta_\x(dz),
\end{equation}
valid for any $f \in C^{sym}(X^2)$ and $\x = (x_1, \ldots, x_n) \in X^n$.

\section{Results}

First we recall some known results on the Cauchy problem for
equation \eqref{eq1.5}. A proof of the following two results can be
found in \cite{No} and \cite{Ko3} respectively. Recall that we
always assume that our continuous coagulation kernel $K(x_1,
x_2;dy)$ preserves $E$ and enjoys the estimate \eqref{eq1.1}.

\begin{proposition}\label{prop2.1}
If a finite measure $\mu_0$ has a finite moment of order $\beta \geq 2$, i.e. if
\begin{equation}\label{eq2.1}
\int_X ( 1 + E^\beta(y)) \mu_0(dy) < \infty,
\end{equation}
then equation \eqref{eq1.5} has a unique solution $\mu_t$ with the
initial condition $\mu_0$ satisfying \eqref{eq1.7} and \eqref{eq1.8} for arbitrary $t$. Moreover,
\begin{equation}\label{eq2.2}
\sup_{s \leq t} \int_X E^\beta(y) \mu_s(dy) \leq c(C,t,\beta, (1 + E, \mu_0))(E^\beta, \mu_0)
\end{equation}
with a constant $c$, and the mapping $\mu_0 \mapsto \mu_t$ is Lipschitz continuous so that
\begin{equation}\label{eq2.3}
\sup_{s \leq t} \| \mu_s(\mu^1_0) - \mu_s(\mu^2_0)\|_{1 + E^\omega}
 \leq c(C,t,\beta,(1+E,\mu_0^1+\mu_0^2))(1 + E^{1 + \omega}, \mu^1_0 + \mu_0^2) \| \mu^1_0 - \mu^2_0\|_{1 + E^\omega}
\end{equation}
for any $\omega \in [1, \beta -1]$.
\end{proposition}

\begin{proposition}\label{prop2.2}
Solutions $\mu_t$ from the previous Proposition enjoy the following
regularity properties:

(i) for any $g \in B_{1 + E^{\beta},\infty}$ (respectively $g \in
B_{1 + E^{\beta -1},\infty}$)
 the function $\int g(x) \mu_t(dx)$ is a continuous function of $t$
 (respectively
continuously differentiable function of $t$ and \eqref{eq1.5}
holds);

(ii) the function $t \mapsto \mu_t$ is absolutely continuous in the norm topology
 of $\MC_{1 + E^{\beta-1}}(X)$ and is continuously differentiable and satisfies
  the strong version of \eqref{eq1.5} in the norm topology of
  $\MC_{1+E^{\beta -\gamma}}(X)$ for any $\gamma \in (1, \beta]$.
\end{proposition}

{\it Remarks.}
\begin{enumerate}
\item The basic ideas of proving Proposition \ref{prop2.1} go back to the analysis of
the Boltzmann equation in \cite{Po}. Formulas \eqref{eq2.2},
\eqref{eq2.3} are proved in \cite{No} only for $\beta =2$ and
$\omega = 1$ respectively, but the above extension is
straightforward.
\item Statement (ii) of Proposition \ref{prop2.1} is proved in \cite{Ko3} only for $\gamma = \beta$,
 but the extension given above is straightforward. In fact (ii)
 is done in the same way as the similar statement of Theorem
\ref{thA2} from Appendix.
\end{enumerate}

It is worth to observe that the operator $L_h$ has the form of the
r.h.s. of equation \eqref{eqA1} from the Appendix with $\MC_{h
\delta}$ instead of $X$ and with the (time homogeneous) intensity
\begin{equation}
\label{eq2.4}
 a(h \delta_\x) = h \sum_{I \subset \{1, \ldots,
n\}:|I|=2} \int K(\x_I;dy) \leq 3Ch^{-1} (1 + E, h \delta_\x)(1,
h\delta_\x).
\end{equation}
As the jumps in \eqref{eq1.4} increase neither $(1, h\delta_\x))$ nor $(E, h \delta_\x)$,
 it is convenient to consider the process $Z^t_h$ on a reduced state space
$$
\MC_{h \delta}^{e_0, e_1} = \{ Y \in \MC_{h \delta} : (1, Y) \leq e_0, (E,Y) \leq e_1\}.
$$
On this reduced space the intensity \eqref{eq2.4} is bounded (not uniformly in $h$).
 Hence $L_h$ is bounded in $C(\MC_{h
\delta}^{e_0, e_1})$ and generates a strongly continuous semigroup of contractions there,
 which we shall denote by $T^h_t$.

Let $T_t$ be a semigroup specified by the solution of \eqref{eq1.5},
 i.e. $T_tf(\mu) = f(\mu_t)$, where $\mu_t$ is the solution of \eqref{eq1.5}
  with the initial condition $\mu$ given by Proposition \ref{prop2.1} with some $\beta \geq 2$.
We can formulate now our first result.

\begin{theorem}\label{th1}[The rate of convergence in LLN]
Let $g$ be a continuous symmetric function on $X^m$ and $F(Y) = (g, Y^{\otimes m})$.
 Assume $Y = h \delta_\x$ belongs to $\MC_{h \delta}^{e_0, e_1}$,
 where $\x = (x_1, \ldots, x_n)$.  Then under the condition (C1) or (C2)
\begin{align}\label{eq2.5}
\sup_{s \leq t} & |T^h_{t}F(Y) - T_tF(Y)| \nonumber \\
& \leq h \kappa (C, m, k, t, e_0, e_1) \|g\|_{(1 + E^k)^{\otimes m}}
(1 + E^{k+3}, Y) (1 + E^3, Y)(1 + E^k, Y)^{m-1}
\end{align}
for any $k\geq 1$ and under the condition (C3)
\begin{align}\label{eq2.6}
\sup_{s \leq t} & |T^h_{t}F(Y) - T_tF(Y)|  \nonumber \\
& \leq h \kappa (C, m, k, t, e_0, e_1) \|g\|_{C^{2,sym}_{(1 +
E^k)^{\otimes m}}(X^m)} (1 + E^{k+4}, Y) (1 + E^{k+1},Y) (1 +
E^k,Y)^{m-1}
\end{align}
for any $k\geq 0$ with a constant $\kappa$.
\end{theorem}

{\it Remarks.}
\begin{enumerate}
\item We give the hierarchy of estimates for the error term making
precise an intuitively clear fact that the power of growth of the
polynomial functions on measures for which LLN can be established
depends on the order of the finite moments of the initial measure.
\item The estimates in case (C2) can be improved. However, not going
 into this detail allows one to keep unified formulae for cases (C1) and (C2).
\item In Section 7 we prove the same estimates \eqref{eq2.5}, \eqref{eq2.6} for more
general functionals $F$ (not necessarily polynomial).
\end{enumerate}

Recall that
$$
F^h_t(Z^h_0, \mu_0) = h^{-1/2} (Z^h_t(Z^h_0) - \mu_t(\mu_0))
$$
is the process of the normalized fluctuations. The main goal of this paper
 is to prove that as $h \rightarrow 0$ this process converges to the generalized
  Gaussian Ornstein-Uhlenbeck (OU) measure-valued process with the (non-homogeneous) generator
\begin{align}\label{eq2.7}
\Lambda_t&F(Y)
= \frac{1}{2} \int \int \int (\delta F(Y), \delta_y - \delta_{z_1} - \delta_{z_2})
 K(z_1, z_2 ; dy) (Y(dz_1)\mu_t(dz_2) + \mu_t(dz_1) Y(dz_2)) \nonumber\\
 & + \frac{1}{4} \int \int \int (\delta^2 F(Y), (\delta_y - \delta_{z_1} - \delta_{z_2})^{\otimes 2})
 K(z_1,z_2; dy)\mu_t(dz_1)\mu_t(dz_2).
 \end{align}

The generalized infinite dimensional Ornstein- Uhlenbeck processes
 and the corresponding Mehler semigroups represent a
 widely discussed topic in the current mathematical literature, see e.g. \cite{LR} and
references therein for general theory, \cite{Ne} for some properties
of Gaussian Mehler semigroups and \cite{Da2} for the connection with
branching processes with immigration. The peculiarity of the process
we are dealing with lies in its 'growing coefficients'. We shall
analyze this process by the analytic tools developed in Sections 4
and 5. Let us start its discussion with an obvious observation that
the polynomial functionals of the form $F(Y) = (g, Y^{\otimes m})$,
$g \in C^{sym}(X^m)$, on measures are invariant under $\Lambda_t$.
In particular, for a linear functional $F(Y) = (g,Y)$
\begin{equation}
\label{eq2.8}
\Lambda_tF(Y) = \frac{1}{2} \int \int \int (g(y) - g(z_1) - g(z_2))
K(z_1,z_2;dy)(Y(dz_1)\mu_t(dz_2) + \mu_t(dz_1) Y(dz_2)).
\end{equation}
Hence the evolution (in the inverse time) of the linear functionals
specified by the equation $\dot F_t = -\Lambda_tF_t$, $F_t(Y) =
(g_t,Y)$ can be described by the equation
\begin{equation}
\label{eq2.9}
\dot{g}(z) = -\Lambda_tg(z) = -\int \int (g(y) - g(x) - g(z))
K(x,z;dy)\mu_t(dx)
\end{equation}
on the coefficient functions $g_t$ (with some abuse of notation we
denoted the action of $\Lambda_t$ on the coefficient functions again
by $\Lambda_t$). Let $U^{t,r}$ be the backward propagator of this
equation, i.e. the resolving operator of the Cauchy problem $\dot{g}
= - \Lambda_tg$ for $t \leq r$ with a given $g_r$.
 As we shall show in Propositions \ref{prop5.2} - \ref{prop5.4}, the
 evolution $U^{t,r}$ is well defined in $C_{1+E^k}(X)$ in cases
 (C1)-(C2), and in $C^{2,0}_{1+E^k}(X)$ in case (C3).

\begin{theorem}\label{th2}[CLT: convergence of linear functionals]
Under condition (C1) or (C2)
\begin{align}
\label{eq2.10}
\sup_{s \leq t} &| \E(g,F^h_s(Z^h_0, \mu_0)) - (U^{0,s}g, F^h_0) | \nonumber \\
& \leq \kappa (C,t,k,e_0, e_1) \sqrt{h} \|g\|_{1 + E^k}(1 + E^{k+3},
Z^h_0 + \mu_0)^2
 \left( 1 + \left\| \frac{Z^h_0 - \mu_0}{\sqrt{h}} \right\|^2_{\MC_{1 + E^{k+1}}(X)}\right)
\end{align}
for all $k \geq 1$, $g\in C_{1+E^k}(X)$, and under condition (C3)
\begin{align}
\label{eq2.11}
\sup_{s \leq t} &| \E(g,F^h_s(Z^h_0, \mu_0)) - (U^{0,s}g, F^h_0) | \nonumber \\
& \leq \kappa (C,t,k,e_0, e_1) \sqrt{h}
\|g\|_{C^{2,0}_{1 + E^k}}(1 + E^{k+5}, Z^h_0 + \mu_0)^3
 \left( 1 + \left\| \frac{Z^h_0 - \mu_0}{\sqrt{h}} \right\|^2_{\MC^1_{1 + E^{k+1}}(X)}\right)
\end{align}
for all $k \geq 0$, $g\in C^{2,0}_{1+E^k}(X)$, where the bald $\E$
denotes the expectation with respect to the process $Z^h_t$.
\end{theorem}

To shorten the exposition, we shall deal in the future only with the
most important case (C3). Though all the results have natural
modifications in cases (C1) and (C2), let us stress again that for
their applicability in cases (C1), (C2) one needs the initial
fluctuation $F_0^h$ to be bounded in the norm of $\MC_{1 +
E^{k+1}}(X)$, which is possible basically only for discrete initial
distributions $\mu_0$.

For our purposes it will be enough to construct the propagator of
the equation $\dot F=-\Lambda_tF$ only on the set of cylinder
functions $\CC_k^n=\CC_k^n(\MC^m_{1+E^k})$, $m=1,2$, on measures
that have the form
 \begin{equation}\label{eq2.12}
\Phi_f^{\phi_1,...,\phi_n}(Y)=f((\phi_1,Y),...,(\phi_n,Y))
\end{equation}
with $f\in C(\R^n)$, and $\phi_1,...,\phi_n \in C^{m,0}_{1+E^k}$. By
$\CC_k$ we shall denote the union of $\CC_k^n$ for all $n=0,1,...$
(of course, functions from $\CC_k^0$ are just constants). Similarly
one defines the cylinder functions $\CC_k^n((L^{m,0}_{2,1+E^k})')$
under condition (C3).

The Banach space of $k$ times continuously differentiable functions
on $\R^d$ (with the norm being the maximum of the sup-norms of a
function and all its partial derivative up to and including the
order $k$) will be denoted, as usual, by $C^k(\R^d)$.

\begin{theorem}\label{th3}[limiting Mehler propagator]
Under the condition (C3) for any $k\ge 0$ and a $\mu_0$ such that
$(1+E^{k+1},\mu_0)<\infty$
 there exists a propagator $OU^{t,r}$ of contractions on ${\CC}_k$ preserving
the subspaces ${\CC}_k^n$, $n=0,1,2,...$ such that $OU^{t,r}F$,
$F\in {\CC}_k$, depends continuously on $t$ in the topology of the
uniform convergence on bounded subsets of $\MC^m_{1+E^k}$, $m=1,2$
(respectively
 $(L^{m,0}_{2,1+E^k})'$ in case $k>1/2$)
and solves the equation $\dot F=-\Lambda_tF$ in the sense that if $f
\in C^2(\R^d)$ in \eqref{eq2.12}, then
\begin{equation}
 \label{eq2.13}
{d \over dt} OU^{t,r}\Phi_f^{\phi_1,..,\phi_n}(Y) =-\Lambda_t
OU^{t,r}\Phi_f^{\phi_1,..,\phi_n}(Y), \quad 0\le t \le r,
\end{equation}
uniformly for $Y$ from bounded subsets of $\MC^m_{1+E^k}$
(respectively
 $(L^{m,0}_{2,1+E^k})'$).
\end{theorem}

Our goal is to prove that this generalized infinite-dimensional
Ornstein-Uhlenbeck (or Mehler) semigroup describes the limiting
Gaussian distributions of the fluctuation process $F_t^h$.

\begin{theorem}\label{th4}[CLT: convergence of semigroups]
 Suppose $k\ge 0$ and $h_0 >0$ are given such that
 \begin{equation}
 \label{eq2.130}
 \sup_{h\le h_0}(1 + E^{k+5}, Z^h_0 +
\mu_0) <\infty.
\end{equation}
(i) Let
 $\Phi \in \CC_k^n(\MC^2_{1+E^k})$ be given by \eqref{eq2.12} with $f\in C^3(\R^n)$
  and all $\phi_j \in C^{2,0}_{1+E^k}(X)$. Then
 \begin{align}
\label{eq2.131}
& \sup_{s \leq t}  |\E\Phi(F_t^h(Z_0^h,\mu_0))-OU^{0,t}\Phi (F_0^h)| \nonumber \\
& \leq \kappa (C,t,k,e_0, e_1) \sqrt{h} \max_j \|
\phi_j\|_{C^{2,0}_{1 + E^k}} \|f\|_{C^3(\R^n)}(1 + E^{k+5}, Z^h_0 +
\mu_0)^3
 \left( 1 + \left\| \frac{Z^h_0 - \mu_0}{\sqrt{h}} \right\|^2_{\MC^1_{1 +
 E^{k+1}}(X)}\right).
\end{align}
(ii) If $\Phi \in \CC_k^n(\MC^2_{1+E^k})$ (with not necessarily
smooth $f$ in the representation \eqref{eq2.12}) and $F_0^h$
converges to some $F_0$ as $h\to 0$ in the $\star$-weak topology of
$\MC^1_{1+E^{k+1}}$, then
\begin{equation}
\label{eq2.14}
\lim_{h\to 0}|\E\Phi(F_t^h(Z_0^h,\mu_0))-OU^{0,t}\Phi (F_0)|=0
\end{equation}
uniformly for $F_0^h$ from a bounded subset of $\MC^1_{1+E^{k+1}}$
and $t$ from a compact interval.
\end{theorem}

\begin{theorem}\label{th5}[CLT: convergence of finite dimensional distributions]
Suppose \eqref{eq2.130} holds, $\phi_1,...,\phi_n \in
C^{2,0}_{1+E^k}(\R_+)$ and $F_0^h \in (L^{2,0}_{2,1+E^{k+2}})'$
converges to some $F_0$ in $(L^{2,0}_{2,1+E^{k+2}})'$, as $h\to 0$.
Then the $\R^n$-valued random variables
$$
\Phi^h_{t_1,...,t_n}=((\phi_1,F_{t_1}^h(Z_0^h,\mu_0)), ...,
(\phi_n,F_{t_n}^h(Z_0^h,\mu_0))), \quad 0<t_1\leq ... \leq t_n,
$$
converge in distribution, as $h\to 0$, to a Gaussian random variable
with the characteristic function
\begin{equation}\label{eq2.15}
g_{t_1,...,t_n}(p_1,...,p_n)=\exp \{i\sum_{j=1}^n
p_j(U^{0,t_j}\phi_j, F_0)-\sum_{j=1}^n
\int_{t_{j-1}}^{t_j}\sum_{l,k=j}^n p_lp_k
 \Pi (s, U^{s,t_l}\phi_l, U^{s,t_k} \phi_k) \, ds \},
\end{equation}
where $t_0=0$ and
\begin{equation}\label{eq2.16}
 \Pi(t,\phi,\psi)=\frac{1}{4} \int \int \int (\phi \otimes \psi ,
 (\delta_y-\delta_{z_1}-\delta_{z_2})^{\otimes 2})
 K(z_1,z_2;dy)\mu_t(dz_1)\mu_t(dz_2).
\end{equation}
In particular, for $t=t_1=...=t_n$ it implies
$$
\lim_{h\to 0}{\bf E} \exp \{i\sum_{j=1}^n (\phi_j, F_t^h)\} =\exp
\{i\sum_{j=1}^n (U^{0,t}\phi_j, F_0)-\sum_{j,k=1}^n \int_0^t
 \Pi (s, U^{s,t}\phi_j, U^{s,t} \phi_k) \, ds \}.
 $$

\end{theorem}

\begin{theorem}\label{th6}[CLT: convergence of the process of fluctuations]
Suppose the conditions of Theorem \ref{th4} hold. (i) For any $\phi
\in C^{2,0}_{1+E^k}(\R_+)$ the real valued processes $(\phi,
F_t^h(Z_0^h,\mu_0))$ converge in the sense of the distribution in
the Skorohod space of c\`adl\`ag functions (equipped with its
standard $J_1$-topology) to the Gaussian process with
finite-dimensional distributions specified by Theorem \ref{th5}.
(ii) The process of fluctuations $F_t^h(Z_0^h,\mu_0)$ converges in
distributions on the Skorohod space of c\`adl\`ag functions
$D([0,T]; (L^{2,0}_{1+E^{k+2}}(\R_+))')$ (with $J_1$-topology),
where
 $(L^{2,0}_{1+E^{k+2}}(\R_+))'$ is considered in its weak topology, to a
Gaussian process with finite-dimensional distributions specified by
Theorem \ref{th5}.
\end{theorem}

\section{Calculations of generators}\label{sec3}

From now on we denote by $\mu_t = \mu_t(\mu_0)$ the solution to
\eqref{eq1.5} given by Proposition \ref{prop2.1} with a $\beta \geq
2$. To begin with, let us extend the action of $T^h_t$ beyond the
space $C(\MC_{h \delta}^{e_0, e_1})$.

\begin{proposition}\label{prop3.1}
For any positive $e_0,e_1$ and $1 \leq l \leq m$ the operator $L_h$ is bounded
 in the space $C_{(1 + E^l,\cdot)^m} (\MC_{h \delta}^{e_0, e_1})$ and defines
  a strongly continuous semigroup there (again denoted by $T^h_t$) such that
\begin{equation}
\label{eq3.1}
\|T^h_t\|_{C_{(1 + E^l,\cdot)^m} (\MC_{h \delta}^{e_0, e_1})} \leq \exp \{ c (C,m,l) e_1t\}.
\end{equation}
\end{proposition}

\noindent {\it Proof.} Let us show that
\begin{equation}\label{eq3.2}
L_h F(Y) \leq c (C,m,l)e_1F(Y)
\end{equation}
for $Y = h \delta_\x$ and $F(Y) = (1 + E^l, Y)^m$. One has
$$
L_hF(Y)
= h \sum_{I \subset \{1, \ldots, n\} : |I|=2} \int
 [(1 + E^l, Y + h(\delta_y - \delta_{\x_I}))^m - (1 + E^l, Y)^m] K(\x_I;dy).
$$
As
\begin{align*}
(1 + E^l, h(\delta_y - \delta_{x_i} - \delta_{x_j})) & \leq h[(E(x_i) + E(x_j))^l - E^l(x_i) - E^l(x_j)] \\
 & \leq hc(l) [E^{l-1}(x_i)E(x_j) + E(x_i)E^{l-1}(x_j)]
 \end{align*}
 and using the obvious inequality $(a+b)^m - a^m \leq c(m)(a^{m-1} b + b^m)$ one obtains
 \begin{align*}
 L_hF(Y)& \leq hc(m,l) \sum_{I \subset \{1, \ldots, n\} : |I|=2}
  [(1 + E^l, Y)^{m-1}  h(E^{l-1}(x_i)E(x_j) + E(x_i) E^{l-1}(x_j)) \\
  & +h^m (E^{l-1}(x_i)E(x_j) + E(x_i)E^{l-1}(x_j))^m] K(\x_I;dy) \\
   & \leq c(C,m,l) \int \int [(1 + E^l, Y)^{m-1}(E^{l-1}(z_1) E(z_2) + E(z_1) E^{l-1}(z_2)) \\
    & + h^{m-1} (E^{l-1}(z_1) E(z_2) + E(z_1) E^{l-1}(z_2))^m](1 + E(z_1) + E(z_2))Y(dz_1)Y(dz_2),
    \end{align*}
    where we used \eqref{eq1.13}. By symmetry it is enough to estimate
    the integral over the set where $E(z_1) \geq E(z_2)$. Consequently $L_hF(Y)$ does not exceed
    \begin{align*}
c \int & [(1 +E^l,Y)^{m-1} E^{l-1}(z_1)E(z_2) + h^{m-1}(E^{l-1}(z_1)E(z_2))^m](1 + E(z_1))Y(dz_1)Y(dz_2) \\
& \leq c(1 + E^l, Y)^m (E,Y) + h^{m-1} c \int E^{m(l-1)+1}(z_1) E^m(z_2)Y(dz_1)Y(dz_2).
\end{align*}
To prove \eqref{eq3.2} it remains to show that the second term in the last expression can be
 estimated by its first term. This follows from the estimates:
\begin{align*}
& (E^m,Y) = h \sum E^m(x_i) \leq h \left( \sum E^l(x_i)\right)^{m/l} = h^{1-m/l}(E^l,Y)^{m/l}, \\
& (E^{m(l-1)+1},Y) \leq h^{-1} (E^{m(l-1)}, Y)(E,Y) \leq h^{-m(1-1/l)}(E^l,Y)^{m(1-1/l)}(E,Y).
\end{align*}
Once \eqref{eq3.2} is proved it follows from \eqref{eq2.4} that
$L_h$ is bounded in
 $C_{(1 + E^l,\cdot)^m} (\MC_{h \delta}^{e_0, e_1})$, and \eqref{eq3.1} follows from Theorem
 \ref{thA1} (or Proposition \ref{propA1}).

The following statement is a straightforward extension of the
previous one.

\begin{proposition}\label{prop3.2}
The statement of Proposition \ref{prop3.1} remains true if instead of the space
$C_{(1 + E^{l},\cdot)^m} $ one takes a more general space
$C_{(1 + E^{l_1},\cdot)^{m_1} \cdots (1 + E^{l_j},\cdot)^{m_j} }$.
\end{proposition}

Next we shall calculate the generator $\LC$ of the deterministic semigroup $T_t$ and compare it with $L_h$.

\begin{proposition}\label{prop3.21}
(i) If $F \in C^1(\MC_{1+E^{\beta},1 + E^{\beta-1}}(X))$, then

\begin{equation}\label{eq3.3}
\frac{d}{dt} T_tF(\mu_0) = \frac{d}{dt} F(\mu_t) = \LC F(\mu_t),
\end{equation}
with
\begin{equation}\label{eq3.4}
\LC F(Y) = \frac{1}{2} \int_X \int_{X \times X} (\delta F(Y;y) -
\delta F(Y;x_1) - \delta F(Y;x_2)) K(x_1,x_2;dy)Y(dx_1)Y(dx_2).
\end{equation}

(ii) If the variational derivative $\delta ^2 F(Y;x,y)$ exists for
$Y\in \MC^+_{1+E^{\beta}}$ and is a continuous function of three
variables ($Y$ taken in its $\star$-weak topology), then for any $Y
= h \delta_\x$
\begin{align}
\label{eq3.5}
L_h&F(Y) - \LC F(Y) = - \frac{h}{2} \int \int (\delta F(Y;y) - 2 \delta F(Y;z))  K(z,z;dy)Y(dz) \nonumber \\
& + h^3 \int_0^1 (1-s) \, ds \sum_{I \subset \{1, \ldots, n\} : |I|
= 2} \int_X (\delta^2 F(Y + sh(\delta_y - \delta_{\x_I}); \cdot ,
\cdot), (\delta_y - \delta_{\x_I})^{\otimes 2})K(\x_I;dy).
\end{align}.

(iii) If $F\in C(\MC (X))$, $Y = h \delta_\x$, then
\begin{align}
\label{eq3.6}
L_hF(Y) &= \frac{1}{2h} \int \int \int [F(Y + h(\delta_y -
\delta_{z_1} - \delta_{z_2})) - F(Y)]
 K(z_1, z_2;dy) Y(dz_1) Y(dz_2) \nonumber \\
& - \frac{1}{2} \int \int [F(Y + h (\delta_y - 2 \delta_z)) - F(Y)]
K(z,z ; dy) Y(dz).
\end{align}
In particular, if $F(Y)=(\phi,Y)$ with a continuous function $\phi$,
then
\begin{align}
\label{eq3.61} L_hF(Y) &= \frac{1}{2} \int \int \int [\phi (y) -
 \phi (z_1) - \phi (z_2)]
 K(z_1, z_2;dy) Y(dz_1) Y(dz_2) \nonumber \\
& - \frac{h}{2} \int \int [\phi (y) - 2 \phi (z)] K(z,z ; dy) Y(dz).
\end{align}
\end{proposition}

\noindent {\it Proof.} (i) Follows from \eqref{eq1.11} and
Proposition \ref{prop2.2}(i).

(ii) Applying \eqref{eq1.10}(a) to \eqref{eq1.4} yields
\begin{align*}
L_h&F(Y) = h^2 \sum_{I \subset \{1, \ldots, n\} : |I| = 2} \int_X
(\delta F(Y; \cdot), \delta_y - \delta_{\x_I})K(\x_I;dy) \nonumber \\
& + h^3 \int_0^1 (1-s) \, ds \sum_{I \subset \{1, \ldots, n\} : |I|
= 2} \int_X (\delta^2 F(Y + sh(\delta_y - \delta_{\x_I}); \cdot ,
\cdot), (\delta_y - \delta_{\x_I})^{\otimes 2})K(\x_I;dy).
\end{align*}
Transforming the first term of the r.h.s. of this equation by
\eqref{eq1.13}, yields \eqref{eq3.5}.

(iii) Is obtained by applying \eqref{eq1.13} directly to
\eqref{eq1.4}.

\begin{proposition}\label{prop3.3}
The backward propagator
$$
U^{h;s,r}_{fl} : C(\Omega^h_r(\MC^{e_0,e_1}_{h \delta})) \mapsto
C(\Omega^h_s (\MC^{e_0,e_1}_{h \delta}))
$$
of the process of fluctuations $F^h_t$ obtained from $Z^h_t$ by the
deterministic linear transformation $\Omega^h_t(Y) = h^{-1/2}(Y -
\mu_t)$, is given by
\begin{equation}
\label{eq3.11}
 U^{h;s,r}_{fl} F = ( \Omega^h_s)^{-1} T^h_{r-s}
\Omega_r^h F,
\end{equation}
where $\Omega_t^h F(Y) = F(\Omega^h_tY)$, and satisfies the equation
\begin{equation}
\label{eq3.10}
 \frac{d}{ds} U^{h;t,s}_{fl} F =
U^{h;t,s}_{fl}\Lambda^h_s F ; \quad t< s < T,
\end{equation}
for $F\in C^3(\MC_{1+E^{\beta},1+E^{\beta-1}}(X))$, where
 \begin{align}
  \label{eq3.9}
  \Lambda^h_tF(Y) &= \Lambda_tF(Y)
  + \frac{\sqrt{h}}{2} \int \int \int (\delta F(Y), \delta_y - \delta_{z_1} - \delta_{z_2})
   K(z_1,z_2;dy)Y(dz_1)Y(dz_2) \nonumber \\
  & - \frac{\sqrt{h}}{2} \int \int (\delta F(Y), \delta_y -  2 \delta_z) K(z,z;dy)(\mu_t + \sqrt{h}Y)(dz) \nonumber \\
  & +  \frac{\sqrt{h}}{4} \int \int \int (\delta^2 F(Y), (\delta_y - \delta_{z_1} - \delta_{z_2})^{\otimes 2})
   K(z_1,z_2;dy)(Y(dz_1)\mu_t(dz_2) + Y(dz_2) \mu_t(dz_1)) \nonumber \\
   & - \frac{h}{4} \int \int \int (\delta^2 F(Y), (\delta_y -2 \delta_z)^{\otimes 2})
    K(z,z;dy)(\mu_t + \sqrt hY)(dz) \nonumber \\
   & + \frac{h}{4} \int \int \int (\delta^2 F(Y), (\delta_y - \delta_{z_1} - \delta_{z_2})^{\otimes 2})
    K(z_1,z_2;dy)Y(dz_1)Y(dz_2) \nonumber \\
   & + \frac{\sqrt h}{4} \int_0^1 (1-s)^2 \, ds \int \int \int
  (\delta^3 F(Y + s\sqrt{h}(\delta_y - \delta_{z_1} - \delta_{z_2}) , \cdot),
  (\delta_y - \delta_{z_1} - \delta_{z_2})^{\otimes 3}) \nonumber \\
   & \qquad \times K(z_1,z_2;dy) (\mu_t+\sqrt h Y)(dz_1)  (\mu_t+\sqrt h
   Y)(dz_2) \nonumber \\
   & - \frac{h^{3/2}}{4} \int_0^1 (1-s)^2 \, ds \int \int
  (\delta^3 F(Y + s\sqrt{h}(\delta_y - 2\delta_z) , \cdot),
  (\delta_y - 2\delta_z)^{\otimes 3})
   K(z,z;dy) (\mu_t+\sqrt h Y)(dz).
  \end{align}
  \end{proposition}

 \noindent {\it Proof.} According to Lemma \ref{lem9.2} the backward propagator
$U^{h;s,r}_{fl}$ is given by \eqref{eq3.11} and satisfies
\eqref{eq3.10} for $F \in C(\Omega_{[0,T]}(\MC^{e_0,e_1}_{h
\delta}))$ (see Lemma \ref{lem9.2} for this notation), where

\begin{equation}
\label{eq3.7}
\Lambda^h_t \psi
= (\Omega^h_t)^{-1} L_h \Omega^h_t \psi - h^{-1/2} \left( \frac{\delta \psi}{\delta Y}, \dot{\mu}_t \right).
\end{equation}

Applying \eqref{eq3.6} yields
\begin{align*}
L_h \Omega^h_t F(Y) &
= \frac{1}{2h} \int \int \int
\left[ F \left( \frac{Y + h(\delta_y - \delta_{z_1} - \delta_{z_2}) - \mu_t}{\sqrt{h}}\right)
 - F \left( \frac{Y - \mu_t}{\sqrt{h}}\right)\right] \\
& \times K(z_1,z_2;;dy) Y(dz_1) Y(dz_2) \\
& - \frac{1}{2} \int \int \left[ F \left( \frac{Y + h(\delta_y - 2 \delta_s) - \mu_t}{\sqrt{h}}\right)
 - F \left( \frac{Y - \mu_t}{\sqrt{h}}\right)\right] K(z,z;dy)Y(dz)
\end{align*}
and consequently
\begin{align}
\label{eq3.8}
(\Omega^h_t)^{-1} & L_h \Omega^h_t F(Y)
 = \frac{1}{2h} \int \int \int \left[ (F(Y + \sqrt{h} (\delta_y - \delta_{z_1} - \delta_{z_2}) - F(Y)\right] \nonumber \\
 & \times K (z_1,z_2;dy)(\sqrt{h}Y + \mu_t) (dz_1) (\sqrt{h} Y + \mu_t)(dz_2) \nonumber \\
  & - \frac{1}{2} \int \int [F(Y+ \sqrt{h} (\delta_y - 2 \delta_z)) - F(Y)] K(z,z;dy)(\sqrt{h} Y + \mu_t)(dz).
  \end{align}
  Applying \eqref{eq1.10} (b) yields
  \begin{align*}
  F&(Y + \sqrt{h}(\delta_y - \delta_{z_1} - \delta_{z_2})) - F(Y)
   = \sqrt{h} (\delta F(Y), \delta_y - \delta_{z_1} - \delta_{z_2})
    + \frac{h}{2} (\delta^2 F(Y), (\delta_y - \delta_{z_1} - \delta_{z_2})^{\otimes 2}) \\
  & + \frac{h^{3/2}}{2} \int_0^1 (1-s)^2
  (\delta^3 F(Y + s\sqrt{h}(\delta_y - \delta_{z_1} - \delta_{z_2})),
  (\delta_y - \delta_{z_1} - \delta_{z_2})^{\otimes 3})\, ds.
  \end{align*}
  Hence developing the r.h.s. of \eqref{eq3.8} in $h$ yields the term at $h^{-1/2}$ of the form
  $$
  \frac{1}{2} \int \int \int (\delta F(Y), \delta_y - \delta_{z_1} - \delta_{z_2})K(z_1, z_2;dy) \mu_t(dz_1) \mu_t(dz_2),
  $$
  the term at $h^0$ being precisely $\Lambda_tF(Y)$ given by \eqref{eq2.7}, plus the remainder terms of order
   at least $h^{1/2}$. As the above term of order $h^{-1/2}$ cancels with the second term
   in \eqref{eq3.7} one obtains \eqref{eq3.9}.

\section{Derivatives with respect to initial data: existence}

This section is devoted to the analysis of the derivatives of the solutions
to equation \eqref{eq1.5} with respect to the initial data. Namely we are
going to study the signed measures defined as
\begin{equation}\label{eq4.1}
\xi_t = \xi_t(\mu_0;x;dz) = \frac{\delta \mu_t}{\delta \mu_0} (\mu_0;x;dz)
 = \lim_{s \rightarrow 0_+} \frac{1}{s}(\mu_t (\mu_0 + s \delta_x) - \mu_t(\mu_0)).
\end{equation}
We will occasionally omit some arguments in $\xi_t$ to shorten the
formulas.

 To motivate the formulation of rigorous results, let us
start with a short formal calculations.
 Differentiating formally equation \eqref{eq1.5} with respect to the initial
 measure $\mu_0$ one obtains for $\xi_t$ the equation
\begin{equation}\label{eq4.2}
\frac{d}{dt} (g, \xi_t) = \int_{X\times X} \int_X (g(y) - g(x_1) - g(x_2)) K(x_1, x_2;dy)\xi_t(dx_1) \mu_t(dx_2).
\end{equation}
Of course, this is by no means a coincidence that this equation is dual to \eqref{eq2.9}.

Introducing the second derivative
\begin{equation}\label{eq4.3}
\eta_t = \eta_t(x,w) = \eta_t(\mu_0;x,w;dz)
 = \lim_{s \rightarrow 0_+} \frac{1}{s} (\xi_t (\mu_0 + s \delta_w;x) - \xi_t(\mu_0;x)),
\end{equation}
and differentiating \eqref{eq4.2} formally one obtains for $\eta_t$ the equation
\begin{align}
\label{eq4.4}
\frac{d}{dt}& (g, \eta_t(x,w;,\cdot)) =\int_{X\times X} \int_X (g(y) - g(x_1) - g(x_2)) K(x_1, x_2;dy) \nonumber \\
& \times [ \eta_t(x,w;dx_1)\mu_t(dx_2) + \xi_t(x;dx_1) \xi_t(w;dx_2)].
\end{align}
The aim of this section is to justify these calculations and to obtain rough estimates for $\xi_t$ and $\eta_t$.

We start our analysis with a result on approximation of the
solutions to kinetic equations by equations
 with bounded kernels. Let us introduce a cut-off kernel $K_n$
 that enjoys the same properties as $K$ and is such that
 $K_n(x_1,x_2;dy) = K(x_1,x_2;dy)$ whenever $E(x_1)+E(x_2) \leq n$ and $K_n (x_1,x_2) \leq Cn$ everywhere.

For convenience, we shall assume $\beta > 3$ everywhere in this section.

\begin{proposition}\label{prop4.1}
Let $\mu_0 \mapsto \mu_t^n$ be the solution, given by Proposition \ref{prop2.1},
 to the equation \eqref{eq1.5} with $K_n$ instead of $K$. Then
  $\mu^n_t \rightarrow \mu_t$ in the norm topology of $\MC_{1 + E^\omega}(X)$
  with $\omega \in [1, \beta -1)$ and $*$-weakly in $\MC_{1+E^\beta}$ uniformly for $t$ from compact sets.
\end{proposition}

{\it Proof.} As the arguments given below use a rather standard trick in the theory of
kinetic equations (similar ideas lead to a proof of Proposition \ref{prop2.1}) we shall give them only for $\omega = 1$.

Let $\sigma^n_t$ denote the sign of the measure $\mu^n_t - \mu_t$
(i.e. the equivalence class of the densities of $\mu^n_t - \mu_t$
with respect to $|\mu^n_t - \mu_t|$ that
 equal $\pm 1$ respectively in positive and negative parts
 of the Hahn decomposition of this measure) so that $| \mu^n_t  - \mu_t| = \sigma^n_t (\mu^n_t - \mu_t)$.
 By Lemma \ref{lem9.3} one can choose a representative of $\sigma^n_t$
 (that we shall again denote by $\sigma^n_t$) in such a way that
\begin{equation}\label{eq4.5}
(1 + E, |\mu^n_t - \mu_t|) = \int_0^t \left( \sigma^n_s (1 + E), \frac{d}{ds} (\mu^n_s - \mu_s)\right)ds.
\end{equation}

Applying \eqref{eq1.5} one obtains from \eqref{eq4.5} that
\begin{align}
\label{eq4.6}
(1 + E, |\mu^n_t - \mu_t|) &
 = \frac{1}{2} \int_0^t ds
  \int [(\sigma^n_s(1 +E))(y) - (\sigma^n_s(1 + E))(x_1) - (\sigma^n_s(1 + E))(x_2)] \nonumber \\
 & \times [K_n(x_1,x_2;dy)\mu^n_s (dx_1) \mu^n_s(dx_2) - K(x_1,x_2;dy)\mu_s(dx_1) \mu_s(dx_2)].
\end{align}
The expression in the last bracket in \eqref{eq4.6} can be rewritten as
\begin{align}\label{eq4.7}
(K_n - K)&(x_1,x_2;dy)\mu^n_s(dx_1)\mu^n_s(dx_2) \nonumber \\
& +K(x_1,x_2;dy)[(\mu_s^n(dx_1) - \mu_s(dx_1)) \mu^n_s(dx_2) +
\mu_s(dx_1)(\mu^n_s(dx_2) - \mu_s(dx_2))].
\end{align}
As $\mu^n_s$ are uniformly bounded in $\MC_{1 + E^\beta}$ and
$$
(1 + E(x_1) + E(x_2)) \int_X (K_n - K) (x_1,x_2;dy) \leq Cn^{- \epsilon}(1 + E(x_1) + E(x_2))^{2 + \epsilon}
$$
for $ 2 + \epsilon \leq \beta$, the contribution of the first term in \eqref{eq4.7} to the r.h.s.
 of \eqref{eq4.6} tends to zero as $n \rightarrow \infty$. The second and the third terms
  in \eqref{eq4.7} ar similar. Let us analyze the second term only. Its contribution to
  the r.h.s; of \eqref{eq4.6} can be written as
\begin{align*}
\frac{1}{2} &\int_0^t ds
 \int [(\sigma^n_s(1 +E))(y) - (\sigma^n_s(1 + E))(x_1) - (\sigma^n_s(1 + E))(x_2)]  \nonumber \\
& \times K(x_1,x_2;dy) \sigma^n_s(x_1) | \mu^n_s (dx_1) - \mu_s(dx_1)| \mu^n_s(dx_2),
\end{align*}
which does not exceed
\begin{align*}
\frac{1}{2} &\int_0^t ds \int [(1 +E)(y) - (1+E)(x_1) + (1+E)(x_2)] \\
& \times K(x_1,x_2;dy) | \mu^n_s (dx_1) - \mu_s(dx_1)| \mu^n_s(dx_2),
\end{align*}
because $(\sigma_s^n(x_1))^2=1$ and $|\sigma_s^n (x_j)|\le 1$,
$j=1,2$. Since $K$ preserves $E$ and \eqref{eq1.1} holds, the latter
expression does not exceed
\begin{align*}
C \int_0^t ds & \int (1 + E(x_2))(1 + E(x_1) + E(x_2)) | \mu^n_s(dx_1) - \mu_s(dx_1)| \mu^n_s(dx_2) \\
& \leq C \int_0^t ds (1 + E, | \mu^n_s - \mu_s|) \| \mu^n_s\|_{1 +
E^2}.
\end{align*}
Consequently by Gronwall's lemma one concludes that
$$
\| \mu^n_t - \mu_t\|_{1+E}
 = (1 + E, | \mu^n_t  - \mu_t|)
 = o(1)_{n \rightarrow \infty} \exp \left\{ t \sup_{s \in [0,t]} \| \mu_s\|_{1 + E^2} \right\}.
 $$
Finally, once the convergence in the norm topology of any
$\MC_{1+E^\gamma}$ with $\gamma > 0$ is established, the $*$-weak
convergence in $\MC_{1+E^\beta}$ follows from the uniform
boundedness of $\mu_n$ and $\mu$ there.

\begin{proposition}\label{prop4.2}

(i) Under the assumptions of Proposition \ref{prop2.1} the backward
propagator $U^{t,r}$ of equation \eqref{eq2.9} is well defined and
is strongly continuous in the space $C_{1+E^{\beta -1}, \infty}(X)$.
Moreover,
 there exists
a unique solution $\xi_t$ to \eqref{eq4.2} in the sense that $\xi_0
= \delta_x$, $\xi_t$ is a $*$-weakly continuous function $\{t\geq
0\} \mapsto \MC_{1 + E^{\beta - 1}}(X)$ and \eqref{eq4.2} holds for
all $g \in C_{1+E}(X)$. Finally,
\begin{equation}\label{eq4.8}
\| \xi_t(\cdot;x)\|_{1 + E^\omega} \leq \kappa (t, \| \mu_0\|_{1 + E^{1+ \omega}})(1 + E^\omega)(x)
\end{equation}
for all $\omega \in [1, \beta -1]$ and some constant $\kappa$,
$\xi_t$ is continuous with respect to $t$ in the norm topology of
$\MC_{1 + E^{\beta - 1 - \epsilon}}$ and is continuously differentiable
 in the norm topology of $\MC_{1+E^{\beta - 2 - \epsilon}}$ for all $\epsilon > 0$.

(ii) If $\xi^n_t$ are defined as $\xi_t$ but from the cut-off
kernels $K_n$, then $\xi^n_t \rightarrow \xi_t$, as $n \rightarrow
\infty$ in the norm topology of $\MC_{1+E^\omega}$ with $\omega \in
[1, \beta - 2)$ and in the $*$-weak topology of $\MC_{1
+E^{\beta-1}}$.

(iii) $\xi_t$ depends Lipschitz continuously on $\mu_0$ in the norm
of $\MC_{1 + E^\omega}$
 for  $\omega \in [1, \beta -2]$ so that
$$
\sup_{s \leq t} \| \xi_s(\mu^1_0) - \xi_s(\mu^2_0)\|_{1 + E^\omega}
 \leq \kappa (C,t,e_0, e_1, (E^{2+\omega}, \mu^1_0 + \mu^2_0))
  \|\mu^1_0 - \mu^2_0\|_{1 + E^{1 + \omega}}(1 + E^{1 + \omega}(x)).
  $$

(iv) $\xi_t$ can be defined by the r.h.s; of \eqref{eq4.1} with the
limit existing in the norm topology of $\MC_{1 + E^\omega}(X)$ with
$\omega \in [1, \beta -1)$ and in the $*$-weak topology of $\MC_{1 + E^{\beta -1}}$.
\end{proposition}

{\it Proof.} (i) Equation \eqref{eq4.2} is dual to \eqref{eq2.9} and
is a particular case of equation \eqref{eqA14} from Appendix with
\begin{equation}\label{eq4.9}
A_tg(x) = \int_X \int_X (g(y) - g(x)) K(z,x;dy) \mu_t(dz),
\end{equation}
and
\begin{equation}\label{eq4.10}
B_tg(x) = \int_X g(z) \int_X K(z,x;dy) \mu_t(dz).
\end{equation}
In the notations of Theorem \ref{thA2} one has in our case
\[
a_t(x)= \int_X \int_X K(z,x;dy) \mu_t(dz) \leq C(1 + E(x)) \|
\mu_t\|_{1+E},
\]
and for all $\omega \leq \beta -1$
\[
\|B_t g\|_{1+E}=\| B_t g/ (1+E)\| \le C \sup_x \{ \frac{\int g(z)
(1+E(x)+E(z))\mu_t (dz)}{1+E(x)} \}
\]
\[
\leq C \|g\|_{1 + E^\omega} \int (1+E^{\omega}(z))(1+E(z))\mu_t (dz)
 \leq 3C\|g\|_{1+E^\omega} \| \mu_t\|_{1
 +E^{\omega + 1}}.
\]

Moreover, as $\omega \ge 1$

\begin{align*}
& A_t(1+E^\omega)(x)
\leq  C \int_X ((E(x) + E(z))^\omega - E^\omega (x))(1 + E(z) + E(x))
\mu_t(dz) \\
& \leq Cc(\omega) \int_X (E^{\omega-1}(x) E(z) + E^\omega(z)) (1 +
E(z) + E(x)) \mu_t(dz) \leq Cc(\omega)(1 + E^\omega)(x) \|\mu_t\|_{1
+ E^{1 + \omega}}.
\end{align*}

Hence the required well-posedness of the dual equations
\eqref{eq2.9} and \eqref{eq4.2} and estimate
 \eqref{eq4.8} for $\omega =\beta -1$ follow from Theorem \ref{thA2} (i), (ii) with $\psi_1=1+
 E^s$, $s \in [1, \beta - 2)$, and $\psi_2 = 1 + E^{\beta -1}$.
  The last statement of (i) follows from Theorem \ref{thA2} (iii).
  Estimate \eqref{eq4.8} for other $\omega \in [1,\beta -1]$ follows
  again from Theorem \ref{thA2} (i) and the estimates for $a_t$ and
  $B_t$ given above.

{\it Remark.} Note that $1+E(x)$ should be of the order $o(1)_{x\to
\infty}(\psi_2/\psi_1)(x)$ in order to fulfill the condition on the
intensity $a_t$ from Theorem \ref{thA1}. Hence the necessity of the
condition $\omega <\beta-2$.

(ii) The proof is the same as the proof of Proposition \ref{prop4.1} above.

(iii) The proof of this statement is practically the same as for the
corresponding statement (see Proposition \ref{prop2.1}(i)) for the
solution of kinetic equation and uses the same trick as in the proof
of Proposition \ref{prop4.1} above. Namely denoting $\xi^j_t =
\xi_t(\mu^j_0), j = 1,2$, one writes
\begin{align*}
\| \xi^1_t - \xi^2_t\|_{1 + e^\omega}&
= \int_0^t ds \int [(\sigma_s(1 + E^\omega))(y) - (\sigma_s(1 + E^\omega))(x_1) - (\sigma_s(1 + E^\omega))(x_2)] \\
& \times K(x_1, x_2;dy)[\xi^1_s(dx_1) \mu^1_s(dx_2) - \xi^2_s(dx_1) \mu^2_s(dx_2)],
\end{align*}
where $\sigma_s$ denotes the sign of the measure $\xi^1_t - \xi^2_t$
 (again chosen according to Lemma \ref{lem9.3}). Next, rewriting
$$
\xi^1_s(dx_1) \mu^1_s(dx_2) - \xi^2_s(dx_1) \mu^2_s(dx_2) = \sigma_s(x_1) |
\xi^1_s - \xi^2_s| (dx_1) \mu^1_s(dx_2) + \xi^2_s(dx_1)(\mu^1_s -
\mu^2_s)(dx_2)
$$
 one estimates from above the contribution of the
first term in the above expression for $\| \xi^1_t - \xi^2_t\|_{1 +
e^\omega}$ by
\begin{align*}
\int_0^t ds &\int[E^\omega(y) - E^\omega(x_1) + E^\omega (x_2) + 1]
 K(x_1,x_2;dy)|\xi^1_s - \xi^2_s| (dx_1) \mu^1_s(dx_2) \\
& \leq c(\omega) C \int_0^t ds \int [E^{\omega-1} (x_1) E(x_2) + E^\omega(x_2) + 1]
(1 + E(x_1) + E(x_2)) |\xi^1_s - \xi^2_s\| (dx_1)\mu^1_s(dx_2) \\
& \leq \kappa  (C, \omega, e_0, e_1) \int_0^t ds \| \xi_s^1 - \xi^2_s\|_{1 + E^\omega} \|\mu^1_s\|_{1 + E^{\omega+1}},
\end{align*}
and the contribution of the second term by
\begin{align*}
\kappa & (C, \omega, e_0, e_1) \int_0^t ds \| \mu^1_s - \mu^2_s\|_{1 + E^{\omega+1}} \| \xi^2_s\|_{1 + E^{\omega+1}} \\
& \leq \kappa (C, \omega, e_0, e_1, (E^{2 + \omega}, \mu^1_0 + \mu^2_0))t
 \| \mu^1_0 -\mu^2_0\|_{1 + E^{\omega +1}} \| \xi^2_0\|_{1 + E^{\omega +1}}.
\end{align*}
It remains to apply Gronwall's lemma to complete the proof of statement (iii).

(iv) General results on the derivatives of the evolution systems
with respect to the initial data seem not to be applied directly for
\eqref{eq1.5}. But they can be applied to the cut-off equations (and
this is the only reason for introducing these cut-offs in our
exposition). Namely, as can be easily seen (this is a simplified
"bounded coefficients" version of Proposition 2.2(ii)), the solution
$\mu^n_t$ to the cut-off version of the kinetic equations
\eqref{eq1.5} satisfies this equation strongly in the norm topology
of $\MC_{1 + E^{\beta-\epsilon}}$ for any $\epsilon > 0$. Moreover,
$\mu^t_n$ depends Lipshtiz continuously on $\mu_0$ in the same
topology, the r.h.s. of the cut-off version of \eqref{eq1.5} is
differentiable with respect to $\mu_t$ in the same topology and
$\xi^n_t$ satisfies the equation in variation \eqref{eq4.2} in the
same topology. Hence it follows from Proposition 6.5.3 of \cite{Mar}
that
$$
\xi^n_t = \xi^n_t(\mu_0;x;dz)
 = \lim_{s \rightarrow 0_+} \frac{1}{s} (\mu^n_t (\mu_0 ,+ s \delta_x) - \mu_t(\mu_0))
$$
in the norm topology of $\MC_{1 + E^{\beta - \epsilon}}$ with $\epsilon > 0$. Consequently
$$
(g, \mu_t^n(\mu_0 + h \delta_x)) - (g, \mu_t^n(\mu_0))
 = \int^h_0  (g, \xi^n_t(\mu_0 + s \delta_x ; x ; \cdot))\, ds
$$
for all $g \in C_{1 + E^{\beta - \epsilon},\infty}(X)$ and $\epsilon > 0$.
Using statement (ii) and the dominated convergence theorem one deduces that
\begin{equation}
\label{eq4.11}
(g, \mu_t(\mu_0 + h \delta_x)) - (g, \mu_t(\mu_0))
 = \int^h_0 (g, \xi_t (\mu_0 + s \delta_x ; x ; \cdot))\, ds
\end{equation}
for all $g \in C_{1+ E^\gamma}(X)$ with $\gamma < \beta - 2$.
Again using the dominated convergence and the fact that $\xi_t$
are bounded in $\MC_{1+E^{\beta-1}}$ (as they are $\star$-weak continuous there)
 one deduces that \eqref{eq4.11} holds for $g\in C_{1 + E^{\beta -1}, \infty}(X)$.
  Next, for these $g$ the expression under the integral in the r.h.s.
  of \eqref{eq4.11} depends continuously on $s$ due to Theorem \ref{thA2} (iv),
   which justifies the weak form of the limit \eqref{eq4.1} (in the $\star$-weak
    topology of $\MC_{1 + E^{\beta -1}}$).
  At last, by statement (iii) $\xi_t$ depends Lipshitz continuously on $s$
  in the r.h.s. of \eqref{eq4.11} in the norm topology of $\MC_{1 + E^\gamma}$
   with $\gamma < \beta - 2$. As $\xi_t$ are bounded in $\MC_{1 + E^{\beta - 2}}$
    it implies that $\xi_t$ depends continuously on $s$ in the r.h.s. of \eqref{eq4.11}
     in the norm topology of $\MC_{1 + E^\gamma}$ with $\gamma < \beta - 1$.
     Hence \eqref{eq4.11} implies \eqref{eq4.1} in the norm topology of $\MC
     _{1+E^{\gamma}}(X)$,
     $\gamma <\beta -1$,
     completing the proof of Proposition \ref{prop4.2}.

\begin{proposition}
\label{prop4.3}
(i) Under the assumptions of Proposition \ref{prop2.1}
there exists a unique solution $\eta_t$ to \eqref{eq4.4}
 in the sense that $\eta_0 = 0$, $\eta_t$ is
  a $*$-weakly continuous function $t \mapsto \MC_{1 + E^{\beta - 2}}$
   and \eqref{eq4.4} holds for $g \in C_{1 + E}(X)$. Moreover
\begin{align}
\label{eq4.12}
\| \eta_t(x,w; \cdot)\|_{1 + E^\omega}& \leq \kappa (C,t, \| \mu_0\|_{1 + E^\beta}) \nonumber \\
& \times \sup_{s \in [0,t]} (\| \xi_s(x;\cdot)\|_{1 + E^{\omega + \alpha}} \| \xi_s(w;\cdot)\|_{1 + E}
 + \| \xi_s(w;\cdot)\|_{1 + E^{\omega + \alpha}} \| \xi_s(x;\cdot)\|_{1 + E})
\end{align}
for $1 \leq \omega \leq \beta -2$ and some $\kappa$.

(ii) If $\eta^n_t$ are defined analogously to $\eta_t$ but
 from the cut-off kernels $K_n$, then $\eta^n_t \rightarrow \eta_t$
  in the norm topology of $\MC_{1 + E^\gamma}$ with $\gamma < \beta - 3$
   and in the $*$-weak topology of $\MC_{1 + E^{\beta-2}}$.

(iii) $\eta_t$ can be defined by the r.h.s. of \eqref{eq4.3} in the
 norm topology of $\MC_{1 + E^\gamma}$ with $\gamma < \beta - 2$
  and in the $*$-weak topology of $\MC_{1 + E^{\beta -2}}$.
\end{proposition}

{\it Proof.} (i) Linear equation \eqref{eq4.4} differs from
 equation \eqref{eq4.2} by an additional non homogeneous term.
  Hence one deduces from Proposition \ref{prop2.1} (i) the
   well posedness of this equation and the explicit formula
\begin{equation}\label{eq4.13}
\eta_t (x,w) = \int_0^t V^{t,s} \Omega_s(x,w) ds,
\end{equation}
where $V^{t,s}$ is a resolving operator to the Cauchy problem of
 equation \eqref{eq4.2} given by Proposition \ref{prop4.2}(i)
  (or directly form Theorem \ref{thA2}) and $\Omega_s(x,w)$ is the measure defined weakly as
\begin{equation}\label{eq4.14}
(g, \Omega_s(x,w)) = \int_{X \times X} \int_X (g(y) - g(x_1) - g(x_2)) K(x_1, x_2;dy) \xi_t(x;dx_1) \xi_t(w;dx_2).
\end{equation}
From this formula and the properties of $\xi_t$ obtained above statement (i) follows.

(ii) This follows from \eqref{eq4.13} and Proposition \ref{prop2.2}(ii).

(iii) As in the proof of Proposition \ref{prop2.2}(iv), we first prove the formula
\begin{equation}\label{eq4.15}
(g, \xi_t(\mu_0 + h \delta_w;x,\cdot)) - (g, \xi_t(\mu_0;x,\cdot)) = \int_0^h (g, \eta_t(\mu_0 + s \delta_w;x,w;\cdot))ds
\end{equation}
for $g \in C_\infty(X)$ by using the approximation $\eta^n_t$, and the dominated
 convergence. Then the validity of \eqref{eq4.15} is extended to all
  $g \in C_{1 + E^{\beta-2},\infty}$ using the dominated convergence and the above obtained
   bounds for $\eta_t$ and $\xi_t$. By continuity of the expression under the integral in the
    r.h.s. of \eqref{eq4.14} we justify the limit \eqref{eq4.3} in the $*$-weak topology
  of $\MC_{1 + E^{\beta-2}}(X)$ completing the proof of Proposition \ref{prop4.3}.

\section{Derivatives with respect to initial data: estimates}\label{sec5}

Straightforward application of Theorem \ref{thA2} of the Appendix would give
 exponential dependence on $(E^\beta, \mu_0)$ of the constant $\kappa$ in \eqref{eq4.8}.
  And this is not sufficient for our purposes. The aim of this Section is to obtain
   more precise estimates for $\xi_t$. Unlike the rough results of the previous section
    that can be more or less straightforwardly extended to very general models with fragmentation,
     collision breakage and their non-binary versions (analyzed in \cite{BK}, \cite{Ko2}, \cite{Ko3}),
  the arguments of this section use more specific properties of the model under consideration.

We shall use the notations of the previous section, assuming in
particular that $A_t$ and $B_t$ are given  by \eqref{eq4.9}, \eqref{eq4.10} respectively.
 Due to the results of the previous section we are able to assume that all the Cauchy
 problems we are dealing with are well-posed. Recall that we denote by $U^{t,r}$ the backward propagator
 of the equation \eqref{eq2.9}.

Let us start with an estimate of the backward propagator $U^{t,r}_A$
 of the equation $\dot{g} = - A_tg$ that holds without additional assumptions (C1)-(C3).

\begin{proposition}\label{prop5.1}
For all $k\geq 0$, $U^{t,r}_A$ is a contraction in
$C_{(1+E^k)^{-1}}$ and
\begin{equation}
 \label{eq5.1}
 |U^{t,r}_Ag(x)| \leq \kappa (C,k, r, e_0, e_1)\|
g\|_{1 + E^k}[(1 + E^k)(x) + (E^{k+1}, \mu_0)].
\end{equation}
\end{proposition}

{\it Proof.} $U_A^{t,r}$ is a contraction in $C_{(1+E^k)^{-1}}$ by
Proposition \ref{propA1}, because $A_t((1+E^k)^{-1}) \le 0$ (and
this holds, because $E^k(y)\ge E^k(x)$ in the support of the measure
$K(z,x;dy)$).
 Next
$$
A_t(1 + E^k)(x) \leq C \int [(E(x) + E(z))^k - E^k(x)](1 + E(x) +
E(z)) \mu_t(dz).
$$
 Using the elementary inequality
$$
((a + b)^k - a^k) (1 + a + b) \leq c(k) (a^k (1 + b) + b^{k+1}+ 1)
$$
that is valid for all positive $a,b,k$ with some constants $c(k)$
yields
$$
A_t(1+ E^k)(x) \leq Cc(k) [E^k(x) (e_0 + e_1) + e_0 + (E^{k+1},
\mu_t)].
$$
 Then by \eqref{eq2.2}
$$
A_t(1+ E^k)(x) \leq \kappa (C,k,t,e_0,e_1) [E^k(x) + 1 + (E^{k+1},
\mu_0)].
$$
Hence \eqref{eq5.1} follows by Lemma \ref{lemA2} and the fact that
$U^{t,r}_A$ is a contraction.

 To simplify formulas we shall often
use the following elementary inequalities :
\begin{align}\label{eq5.2}
(a) \hspace{1cm} &(E^l, \nu)(E^k, \nu) \leq 2(E^{k+l-1}, \nu)(E, \nu),  \nonumber \\
(b) \hspace{1cm}&(E^k, \nu)E(x) \leq (E^{k+1}, \nu) + (E,\nu)E^k(x)
\end{align}
valid for arbitrary positive $\nu$ and $k,l \geq 1$.

\begin{proposition}\label{prop5.2}
Under condition (C1) suppose $k \geq 1$. Then

\begin{equation}\label{eq5.3}
|U^{t,r}g(x)| \leq \kappa (C,k,r, e_0, e_1) \|g\|_{1 + E^k} [1 +
E^k(x) + (E^{k+1}, \mu_0)(1 + E(x))],
\end{equation}
\begin{equation}\label{eq5.4}
\sup_{s \leq t} \| \xi_s(\mu_0;x,\cdot)\|_{1 + E^k} \leq \kappa
(C,t,e_0, e_1)[1 + E^k(x) + (1 + E(x))(E^{k+1}, \mu_0)],
\end{equation}
and
\begin{align}
\label{eq5.5}
\sup_{s \leq t} &\| \eta_s(\mu_0;x,w;\cdot)\|_{1 + E^k} \leq \kappa (C,k,t,e_0,e_1) \nonumber
\\
& \times [(1 + E^{k+1}(x) + (E^{k+1}, \mu_0)(1 + E^2(x)) + (E^{k+3}, \mu_0)(1 + E(x)))(1 + E(w)) \nonumber \\
& + (1 + E^{k+1}(w) + (E^{k+1}, \mu_0)(1 + E^2(w)) + (E^{k+3}, \mu_0)(1 + E(w)))(1 + E(x))].
\end{align}
\end{proposition}

{\it Proof.} The simplicity of condition (C1) stems from the observation that
 the two dimensional functional space generated by the function $E$ and constants
  is invariant under both $A_t$ and $B_t$, and also the full image of $B_t$ belongs to
   this space. Hence representing the solution to $\dot{g} = - (A_t - B_t)g$ as
\begin{equation}\label{eq5.6}
g = U^{t,r}_A g_r + \tilde{g}
\end{equation}
one finds that $\tilde{g}$ belongs to the above mentioned two dimensional space and satisfies the equation
\begin{equation}\label{eq5.7}
\dot{\tilde{g}} = - (A_t - B_t) \tilde{g} + B_t U^{t,r}_A g_r, \qquad \tilde{g}\mid_{t=r} = 0, \qquad t \leq r.
\end{equation}
The corresponding homogeneous Cauchy problem
$$
\dot{\phi} = -(A_t - B_t) \phi, \qquad \phi_r = \alpha + \beta E,
$$
can be written as
$$
\dot \alpha_t +\dot \beta_t E(x)=C\alpha_t(e_1+(1, \mu_t) E(x)),
\quad \alpha_r=\alpha, \beta_r=\beta
$$
in terms of $\phi =\alpha_t +\beta_t E(x)$ and clearly solves
explicitly as
$$
\phi_t = \alpha e^{-e_1(r-t)} + \left[ \beta + \alpha \int_t^r (1, \mu_s)e^{-e_1(r-s)}ds \right] E(x),
$$
which implies that
$$
\| \phi_t\|_{1 + E} \leq  \kappa (r,e_0) \| \phi_r\|_{1 + E}.
$$
It follows from \eqref{eq5.1} that
\begin{align}\label{eq5.8}
|B_t U^{t,r}_A g_r(x)| &\leq  \kappa (C,r, e_0, e_1) \| g_r\|_{1 +
E^k} [B_t (1 + E^k)+(E^{k+1}, \mu_0)B_t1](x)
 \nonumber \\
& \leq \kappa (C, r, e_0, e_1)
 \| g_r\|_{1 + E^k}(1 + (E^{k+1}, \mu_0))(1 + E(x)).
\end{align}
Solving the non-homogeneous equation \eqref{eq5.7} by the Du Hamel principle and using
 the representation \eqref{eq5.6} yields \eqref{eq5.3}. But by duality one gets
\[
 \| \xi_s(\mu_0;x,\cdot)\|_{1 + E^k}=\sup \{(g,\xi_s(\mu_0;x,\cdot)): \|g\|_{1 +
 E^k} \le 1\}
 \]
 \[
 =\sup \{ (U^{0,s}g, \delta_x):\|g\|_{1 +
 E^k} \le 1\}=\sup \{ U^{0,s}g(x):\|g\|_{1 +
 E^k} \le 1\},
 \]
which implies \eqref{eq5.4}.

Now from \eqref{eq4.13}
\begin{align*}
\sup_{s \leq t} & \| \eta(\mu_0;x,w;\cdot)\|_{1 + E^k} \leq t \sup_{s \leq t} \| V^{t,s} \Omega_s(x,w)\|_{1 + E^k}
 = t \sup_{s \leq t} \sup_{|g| \leq 1 + E^k} (U^{s,t}g, \Omega_s(x,w)) \\
& \leq \kappa (C,t,e_0,e_1) \sup_{s \leq t} \sup \{(g, \Omega_s(x,w)): |g(y)| \leq 1 + E^k(y)
 +(1+E(y))(E^{k+1}, \mu_0)\} \\
& \leq \kappa (C,k,t,e_0,e_1) \sup_{s \leq t} \int \int [1 + E^k(x_1) + E^k(x_2)
 + (E^{k+1}, \mu_0)(1 + E(x_1) + E(x_2))] \\
 & \qquad (1 + E(x_1) + E(x_2)) \xi_s(x;dx_1) \xi_s(w;dx_2).
\end{align*}
Dividing this integral into two parts with $E(x_1) \geq E(x_2)$ and $E(x_1) \leq E(x_2)$
 one can estimate the first part as

\begin{align*}
 & \kappa \sup_{s \leq t} \int \int [1 + E^k(x_1) + (E^{k+1}, \mu_0)(1 + E(x_1))](1 + E(x_1))
  \xi_s(x;dx_1)\xi_s(w;dx_2) \\
  & \leq \kappa \sup_{s \leq t}  \|\xi_s(w;\cdot)\|
\left(\|\xi_s(x;\cdot)\|_{1 + E^{k+1} }
 +(E^{k+1}, \mu_0) \| \xi_s(x;\cdot)\|_{1 + E^2}\right) \\
  &\leq \kappa (1 + E(w))[1 + E^{k+1}(x)
+ (1 + E(x))(E^{k+2},\mu_0) \\
   & \qquad + (E^{k+1}, \mu_0)(1 + E^2(x) + (1 + E(x))(E^3, \mu_0)] \\
   &\leq \kappa (1 + E(w))[1 + E^{k+1}(x) + (E^{k+1}, \mu_0)(1 + E^2(x)) + (E^{k+3}, \mu_0)(1 + E(x))],
   \end{align*}
   where we used both \eqref{eq5.2}(a) and \eqref{eq5.2}(b).
    As the integral over the second part is estimated similarly one arrives at \eqref{eq5.5}.

   \begin{proposition}\label{prop5.3}
   Under condition (C2)
   \begin{equation}\label{eq5.9}
   \|U^{t,r}\|_{C_{1 + \sqrt{E}}} \leq \exp \{ 4 C(t-r)(e_0 + e_1)\}
   \end{equation}
   and the estimates \eqref{eq5.3}-\eqref{eq5.5} hold for all $k \geq 1$.
   \end{proposition}

   {\it Proof.}
   Since
   \begin{align*}
   A_t(1 + \sqrt{E})(z) & = \int \int (\sqrt{E(z) + E(x)} - \sqrt{E(z)}) K(z,x;dy) \mu_t (dx) \\
    & \qquad \leq C \int_X \sqrt{E(x)} (1 + \sqrt{E(z)}) (1 + \sqrt{E(x)}) \mu_t (dx) \\
     & \leq C (1 + \sqrt{E(z)}) (\sqrt{E} + E, \mu_t) \leq C(e_0 + 2e_1)(1 + \sqrt{E(z)}),
     \end{align*}
     according to Proposition \ref{propA1} the positivity preserving backward
      propagator $U^{r,t}_A$ of the equation $\dot{g} = - A_tg$ is bounded
       in $C_{1 + \sqrt{E}}(X)$ with the norm not exceeding $\exp \{ C(t-r)(e_0 + 2e_1)\}$. On the other hand
     \begin{align*}
     B_t(1 + \sqrt{E})(z) & \leq C \int (1 + \sqrt{E(x)})^2 (1 + \sqrt{E(z)}) \mu_t (dx) \\
      & \leq 2 C(e_0 + e_1)(1 + \sqrt{E(z)}).
      \end{align*}
Hence $B_t$ are uniformly bonded in $C_{1 + \sqrt{E}}(X)$ with the
norm not exceeding $2C(e_0 + e_1)$. Hence \eqref{eq5.9} follows from
the series representation \eqref{eqA16} for the backward propagator
$U^{r,t}$ of the equation $\dot{g} = - (A_t - B_t)g$.

      Now we use the same arguments as in the proof of Proposition \ref{prop5.2} with
 $$
|B_t U^{t,r}_A g_r(x)| \leq \kappa (C,t - r, e_0, e_1) \|g_r\|_{1+ E^k}(1 + \sqrt{E(x)}) (1 + (E^{k+1}, \mu_0))
 $$
instead of \eqref{eq5.8}. Namely, in the representation of the
solutions to $\dot{g} = - (A_t - B_t)g$ by the series \eqref{eqA16}
the first term is independent of $B_t$ and all other terms belong to
$C_{1 + \sqrt{E}}(X)$ and applying the above estimates for
$U_A^{t,r}$ and $B_t$ in this space one deduces \eqref{eq5.3}. Other
estimate follows now straightforwardly as in the previous
Proposition (even with some improvements that we do not take into
account).

\begin{proposition}\label{prop5.4}
Under condition (C3) for any $k \geq 0$ the spaces $C^{1,0}_{1+E^k}$
and $C^{2,0}_{1+E^k}$ (see Introduction for these notations) are
invariant under $U^{t,r}$ and
\begin{align}\label{eq5.10}
(a) \hspace{1cm} & |(U^{t,r}g)'(x)| \leq \kappa (C,r, k, e_0, e_1)
 \|g\|_{C^{1,0}_{1+E^k}} (1 + E^k(x) + (E^{k+1}, \mu_0)), \nonumber \\
(b) \hspace{1cm} & |(U^{t,r}g)''(x)| \leq \kappa (C,r, k, e_0, e_1)
 \|g\|_{C^{2,0}_{1+E^k}} (1 + E^k(x) + (E^{k+1},
\mu_0)),
\end{align}
\begin{equation}\label{eq5.11}
\sup_{s \leq t} \| \xi_s(\mu_0;x;\cdot)\|_{\MC^1_{1 + E^k}}
 \leq \kappa (C,r, k, e_0, e_1)[E(x) (1 + (E^{k+1}, \mu_0)) + E^{k+1}(x)],
\end{equation}
and
\begin{align}\label{eq5.12}
\sup_{s \leq t}& \| \eta_s(\mu_0;x,w;,\cdot)\|_{\MC^2_{1 + E^k}}
 \leq \kappa (C,t , k, e_0, e_1) (1 + (E^{k+1}, \mu_0)) \nonumber \\
 & \times [(E(x)(1 + E^{k+2}, \mu_0) + E^{k+2}(x))E(w) + (E(w)(1 + E^{k+2}, \mu_0) + E^{k+2}(w))E(x)].
 \end{align}
\end{proposition}

{\it Proof.} Notice first that if $g_r(0) = 0$, then $g_t = 0$ for
all $t$ according to the evolution described by the equation
 $\dot{g} = - (A_t - B_t)g$. Hence the space of functions vanishing
  at the origin is invariant under this evolution.

  Recall that $E(x)=x$ in case (C3). Differentiating
   the equation $\dot{g} = - (A_t - B_t)g$ with respect to the space variable $x$ leads to the equation
\begin{equation}\label{eq5.13}
\dot{g}'(x) = - A_t(g')(x) - \int (g(x + z) - g(x) - g(z)) \frac{\partial K}{\partial x}(x,z) \mu_t(dz).
\end{equation}
For functions $g$ vanishing at the origin this can be rewritten as
$$\dot{g}'(x) = - A_tg' - D_tg'$$
with
$$
D_t \phi(x) = \int \left( \int_x^{x+z} \phi(y)dy - \int_0^z
\phi(y)dy \right) \frac{\partial K}{\partial x}(x,z) \mu_t(dz).
$$
Since
$$
 \|D_t \phi\| \leq 2C \| \phi \| (E , \mu_t)
= 2 Ce_1 \| \phi\|,
 $$
and $U^{t,r}_A$ is a contraction, it follows from representation \eqref{eqA16} with $D_t$ instead of $B_t$ that
$$
\| U^{t,r}\|_{C^{1,0}_1(X)} \leq \kappa (C, r - t, e_0, e_1),
$$
proving \eqref{eq5.10}(a) for $k = 0$. Next, for $k >0$
\begin{align*}
|D_t \phi(x)| & \leq C \| \phi\|_{1 + E^k} \int ((x + z)^{k+1} - x^{k+1} + z^{k+1} + 2z) \mu_t (dz) \\
& \leq C c(k) \| \phi\|_{1 + E^k} \int (x^k z + z^{k+1} + 2 z) \mu_t
(dz),
\end{align*}
which by \eqref{eq2.2} does not exceed
$$
 c(C,k,e_1) \| \phi\|_{1 + E^k}[(1 + x^k) + (E^{k+1}, \mu_0)].
 $$
Hence by Proposition \ref{prop5.1}
$$
\int_t^r | U^{t,s}_A D_s U^{s,r}_Ag(x)|ds
 \leq (r-t) \kappa (C,r,k,e_0,e_1) \|g\|_{1 + E^k} [1 + E^k(x) + (E^{k+1}, \mu_0)],
 $$
which by induction implies
\begin{align*}
\int_{t \leq s_1 \leq \cdots \leq s_n \leq r} &| U^{t,s_1}_A D_{s_1} \cdots D_{s_n} U^{s_n,r}_A g(x) | ds_1 \cdots ds_n
\\
& \leq \frac{(r-t)^n}{n!} \kappa^n (C,r,k,e_0,e_1) \|g\|_{1+E^k} [1 + E^k(x) + (E^{k+1}, \mu_0)].
\end{align*}
Hence \eqref{eq5.10}(a) follows from the representation \eqref{eqA16} to the solution of \eqref{eq5.13}.

Differentiating \eqref{eq5.13} leads to the equation
\begin{equation}\label{eq5.14}
\dot{g}''(x) = - A_t (g'')(x) - \psi_t,
\end{equation}
where
\begin{align*}
\psi_t & = 2 \int (g'(x+z) - g'(x)) \frac{\partial K}{\partial x}(x,z) \mu_t(dz) \\
 & +\int \left( \int_x^{x+z} g'(y) dy - \int_0^z g'(y) dy \right)
  \frac{\partial^2K}{\partial x^2} (x,z) \mu_t(dz).
 \end{align*}
 We know already that for $g_r \in C^2_{1 + E^k}$ the function
  $g'$ belongs to ${1 + E^k}$ with the bound given by \eqref{eq5.10}(a).
  Hence by the Du Hamel principle the solution to \eqref{eq5.14} can be represented as
 $$g''_t = U^{t,r}_A g''_r + \int_t^r U^{t,s}_A \psi_s ds.$$
 As
 $$| \psi_t(x)| \leq \kappa (C,r-t, e_0, e_1)(1 + E^k(x) + (E^{k+1}, \mu_0)),$$
 \eqref{eq5.10}(b) follows, completing the proof of \eqref{eq5.10}, which by duality implies \eqref{eq5.11}.

Next, arguing as in the proof of Proposition \ref{prop5.2} one gets
$$
\sup_{s\leq t} \|\eta_s(\mu_0;x,w;.)\|_{\MC^2_{1+E^k}}
 \leq t\sup_{s\leq t}
 \sup \{ |(U^{s,t}g,\Omega_s(x,w))|: \|g\|_{C^{2,0}_{1+E^k}} \leq 1\}
$$
$$
\leq \kappa (C,t,e_0,e_1) \sup_{s\leq t} \sup_{g\in \Pi_k} (g,\Omega_s(x,w)),
$$
where
$$
\Pi_k=\{ g: g(0)=0, \max(|g'(y)|,|g''(y)|)\leq 1+E^k(y)+(E^{k+1},\mu_0)\}.
$$
It is convenient to introduce a two times continuously differentiable function
$\chi$ on $\R$ such that $\chi(x) \in [0,1]$ for all $x$, and $\chi(x)$ equals
one or zero respectively for $x \geq 1$ and $x \leq -1$. Then write
 $\Omega_s = \Omega^1_s + \Omega^2_s$ with $\Omega^1$ (respectively $\Omega^2$)
  being obtained by \eqref{eq4.14} with $\chi(x_1 - x_2) K(x_1,x_2)$
  (respectively $(1 - \chi(x_1 - x_2))K(x_1,x_2))$ instead of $K(x_1,x_2)$.
If $g\in \Pi_k$, one has
 $$
(g, \Omega^1 _s(x,w)) = \int \int (g(x_1+x_2)-g(x_1)-g(x_2)) \chi (x_1 - x_2) K(x_1, x_2) \xi_s(x;dx_1)\xi_s(w;dx_2),
$$
 which is bounded in magnitude by
 \begin{align*}
 \|\xi_s(w,\cdot)\|_{\MC^1_1(X)} & \sup_{x_2}\left| \frac{\partial}{\partial x_2} \int [(g(x_1+x_2)-g(x_1)-g(x_2))
 \chi(x_1-x_2) K(x_1,x_2)] \xi_s (x;dx_1) \right| \\
& \leq  \|\xi_s(w,\cdot)\|_{\MC^1_1(X)}
  \|\xi_s(x,\cdot)\|_{\MC^1_{1+E^{k+1}}(X)} \\
 & \times \sup_{x_1,x_2} \left| (1 + E^{k+1}(x_1))^{-1}
  \frac{\partial^2}{\partial x_2 \partial x_1} [(g(x_1 + x_2) - g(x_1) - g(x_2)) \chi(x_1 - x_2)K(x_1,x_2)]\right|.
 \end{align*}
Since
\begin{align*}
\frac{\partial^2}{\partial x_2 \partial x_1}&[(g(x_1 + x_2) - g(x_1) - g(x_2)) \chi(x_1 - x_2)K(x_1,x_2)] \\
& = g''(x_1 + x_2) (\chi K)(x_1,x_2)+(g'(x_1+x_2)-g'(x_2)) \frac{\partial (\chi K)(x_1,x_2)}{\partial x_1} \\
& + (g'(x_1 + x_2) - g'(x_1)) \frac{\partial (\chi K)(x_1,x_2)}{\partial x_2}
 + (g(x_1 + x_2) - g(x_1) - g(x_2)) \frac{\partial^2 (\chi K)(x_1,x_2)}{\partial x_1 \partial x_2},
\end{align*}
this expression does not exceed in magnitude
 $C(1+E^{k+1}(x_1)+(E^{k+1},\mu_0)(1+E(x_1))$
 (up to a constant multiplier). Consequently
$$
|(g, \Omega^1_s(x,w)| \leq \kappa (C)
\| \xi_t(w,\cdot)\|_{\MC^1_1(X)} \| \xi_t(x,\cdot)\|_{\MC^1_{1 + E^{k+1}}(X)}(1 + (E^{k+1}, \mu_0)).
$$
Of course, the norm of $\Omega^2_s$ is estimated in the same way.
 Consequently \eqref{eq5.11} leads to \eqref{eq5.12} and completes the proof of Proposition \ref{prop5.4}.

 We shall prove now the Lipschitz continuity of the solutions
  to our kinetic equation with respect to initial
 data in the norm-topology of the space $\MC^1_{1+E^k}$.

\begin{proposition}\label{prop5.5}
Under the condition (C3) for $k\geq 0$ and $m=1,2$
\begin{equation}\label{eq5.15}
\sup_{s \leq t} \| \mu_s(\mu^1_0) - \mu_s(\mu^2_0)\|_{\MC^m_{1 +
E^k}}
 \leq \kappa (C,t,k,e_0,e_1)(1+E^{1+k},\mu^1_0+\mu_0^2)\|\mu^1_0-\mu^2_0\|_{\MC^m_{1 + E^k}}
\end{equation}
\end{proposition}

{\it Proof.} By \eqref{eq4.1} and \eqref{eq1.9}
\begin{equation}\label{eq5.16}
(g,\mu_t(\mu_0^1)-\mu_t(\mu_0^2))
 =\int_0^t ds \int\int g(y)\xi_t(\mu_0^2+s(\mu_0^1-\mu_0^2);x;dy)(\mu_0^1-\mu_0^2)(dx).
 \end{equation}
 Since
 $$
 (g,\xi_t(Y;x;.))=(U^{0,t}g,\xi_0(Y,x;.))=(U^{0,t}g)(x),
 $$
 it follows from Proposition \ref{prop5.4} that $(g,\xi_t(Y;x;.))$ belongs to
  $C^{m,0}_{1+E^k}$ as a function of $x$ whenever $g$ belongs to this space and that
  $$
 \| (g,\xi_t(Y;x;.)) \|_{C^{m,0}_{1+E^k}(X)}
  \leq \kappa (C,t,k,e_0,e_1)\| g\|_{C^{m,0}_{1+E^k}(X)}(1+(E^{k+1},Y)).
  $$
  Consequently \eqref{eq5.15} follows from \eqref{eq5.16}.

We shall discuss now the $L^2$-version of our estimates.

\begin{proposition}\label{prop5.6} Under condition (C3) assume
 $f$ is a positive either non-decreasing or bounded function. Then
$U_A^{t,r}$ are contractions in $L_{2,1/f}$. (Thus $U_A^{t,r}$ yield
natural examples of sub-Markovian propagators with growing
coefficients.)
\end{proposition}

{\it Proof.} First observe that
\begin{equation}\label{eq5.17}
 \int_0^{\infty} (u(x+y)-u(x))g^2(x)u(x)dx \le 0
 \end{equation}
 for any $y\ge 0$ and a non-decreasing non-negative $g$ (and any $u$, if only the integral is
 well defined). In fact, it is equivalent to
 $$
 (T_yu,u)_{L_{2,1/g}} \le (u,u)_{L_{2,1/g}},
 $$
 where $T_yu(x)=u(x+y)$, which in turn follows (by Cauchy
 inequality) from $(T_yu,T_yu)_{L_{2,1/g}}\le (u,u)_{L_{2,1/g}}$.
 And the latter holds, because
 $$
 (T_yu,T_yu)_{L_{2,1/g}}=\int_0^{\infty} u^2(x+y) g^2(x)\, dx
 $$
 $$
 =\int_y^{\infty} u^2(z) g^2(z-y)\, dz \le \int_y^{\infty} u^2(x) g^2(x)\,
 dx.
 $$
 Assume now that $f$ is non-decreasing (and positive). From \eqref{eq5.17} it follows that
for arbitrary $y>0$
\begin{equation}
\label{eq5.18}
 \int_0^{\infty} (u(x+y)-u(x))K(x,y)f^2(x)u(x)\, dx \le 0
 \end{equation}
 (we used here the assumed monotonicity of the kernel $K$),
 and hence $(A_tu,u)_{L_{2,1/f}} \le 0$. Hence $U_A^{t,r}$ can not
 increase the norm of $L_{2,1/f}$. To conclude that it is actually a
 semigroup of contractions it remains to observe that due to
 Proposition
 \ref{prop5.1} there exists a dense subspace in $L_{2,1/f}$ that is
 invariant under $U_A^{t,r}$. Assume now that $f$ is bounded. We
 again have to show the validity of \eqref{eq5.18} for a dense
 invariant subspace of functions $u$. First note that as the
 evolution $U_A^{t,r}$ is well defined on continuous functions and
 preserves positivity and differentiability (by Propositions \ref{prop5.1}
 and \ref{prop5.4})
 it is suffice to show \eqref{eq5.18} for positive functions $u$
 with bounded variation. Assuming that this is the case one can
 represent positive $u$ as the difference $u=u^+-u^-$ of two
 positive non-decreasing functions (by decomposing its derivative
 in its positive and negative parts). As $-(u^-(x+y)-u^-(x)) \le 0$,
 to show \eqref{eq5.18} one needs to show that
 $$
 \int_0^{\infty} (u^+(x+y)-u^+(x))K(x,y)f^2(x)u(x)\, dx \le 0,
 $$
 and as $-u^-$ is negative this in turns follows from
$$
 I_y=\int_0^{\infty} (u^+(x+y)-u^+(x))K(x,y)f^2(x)u^+(x)\, dx \le 0.
 $$
 Denoting by $M$ an upper bound for $f^2$ we can write
 $$
 I_y=\int_0^{\infty} (u^+(x+y)-u^+(x))K(x,y)Mu(x)\, dx
  -\int_0^{\infty} (u^+(x+y)-u^+(x))K(x,y)(M-f^2(x))u(x)\, dx,
 $$
 which is negative, because the integrand in the second term is
 positive and the first term is negative by \eqref{eq5.17}.

We shall consider now equation (5.13) that can be written in the
form
\begin{equation}\label{eq5.19}
\dot{g}'(x) = - A_t(g')(x) - D_t^1g'-D_t^2g',
\end{equation}
where
$$
(D_t^1\phi)(x)=\int_0^x \left( \int_x^{x+z}\phi (y) \, dy \frac
{\partial K} {\partial x} (x,z) \right) \mu_t (dz),
$$
$$
(D_t^2\phi)(x)=\int_0^{\infty} \left( {\bf 1}_{x<z}\int_x^{x+z}\phi
(y) \, dy -\int_0^z \phi (y) \, dy \right)\frac {\partial K}
{\partial x} (x,z) \mu_t (dz).
$$

\begin{proposition}\label{prop5.7}
Under condition (C3) for any $f_k(x)=1+x^k$, $k>1/2$, the spaces
$L^{m,0}_{2,f_k}$, $m=1,2$ (see introduction for these notations),
are invariant under $U^{t,r}$ and
 $$
 \|U^{t,r}\|_{L^{m,0}_{2,f_k}}\le \kappa (C,r,k,e_0,e_1)
 (1+(E^{k+1/2},\mu_0)), \quad m=1,2.
 $$
 Moreover, for $g\in L^{m,0}_{2,f_k}$ one can represent
 $(U^{t,r}g)'$ as the sum of a function from $L^{m,0}_{2,f_k}$ with
 the norm not exceeding $\kappa (C,r,k,e_0,e_1)\|g\|_{L^{m,0}_{2,f_k}}$ and a uniformly
 bounded function with the sup-norm not exceeding
 $\kappa
 (C,r,k,e_0,e_1)\|g\|_{L^{m,0}_{2,f_k}}(1+(E^{k+1/2},\mu_0))$.
 Consequently
 \begin{equation}\label{eq5.20}
\sup_{s \leq t} \| \xi_s(\mu_0;x;\cdot)\|_{(L^{1,0}_{2,f_k})'}
 \leq \kappa (C,r, k, e_0, e_1)[E(x)(E^{k+1}, \mu_0) +(1+
 E^{k+1/2}(x)].
\end{equation}

\end{proposition}

{\it Proof.} Let us show first that
\begin{equation}
\label{eq5.21}
 \|D_t^1\|_{L_{2,f_k}} \le e_1 2^kC.
 \end{equation}
 In fact,
 for a continuous positive $\phi$ and an arbitrary $z>0$
 $$
 \|{\bf 1}_{z\le x} \int_x^{x+z} \phi (y) \, dy \|^2_{L_{2,f_k}}
 = \lim_{n\to \infty} \int {\bf 1}_{z\le x}\sum_{i,j=1}^n \phi
 (x+{jz \over n})\phi (x+{iz \over n})\frac {z^2}{n^2} f_k^{-2}(x) \,dx
 $$
 $$
 \le z^2\lim_{n\to \infty} \int {\bf 1}_{z\le x}\sum_{i,j=1}^n
   (\phi/f_k) (x+{jz \over n})(\phi/f_k) (x+{iz \over n})
   \frac {1}{n^2} 2^{2k} \, dx,
   $$
because
$$
\frac{1}{f_k(x)} \le \frac{2^k}{f_k(2x)}\le \frac{2^k}{f_k(x+jz/n)}
$$
for all $j\le n, z\le x$. Taking now into account that
$$
\int {\bf 1}_{z\le x}\sum_{i,j=1}^n
   (\phi/f_k) (x+{jz \over n})(\phi/f_k) (x+{iz \over n})\, dx
   \le \|\phi/f_k\|^2_{L_2}
   $$
   one deduces that
$$
 \|{\bf 1}_{z\le x} \int_x^{x+z} \phi (y) \, dy \|^2_{L_{2,f_k}}
\le z^2 2^{2k} \|\phi \|^2_{L_{2,f_k}}.
   $$
   Consequently
\[
\|D_t^1 \phi \|_{L_{2,f_k}} \le C \int_0^{\infty} \|{\bf 1}_{z\le x}
\int_x^{x+z} \phi (y) \, dy \|_{L_{2,f_k}} \mu_t (dz) \le C 2^k e_1
\| \phi \|_{L_{2,f_k}},
\]
which implies \eqref{eq5.21}.

   As clearly the same bounds hold for
 $\|D_t^1\|_{C(\R_+)}$, the equation
 $$
 \dot{g}'(x) = - A_t(g')(x) - D_t^1g'
 $$
 specifies a propagator $\tilde U^{t,r}$, $t\in [0,r]$, of bounded
 operators in both $C({\R}_+)$ and $L_{2,f_k}(\R_+)$ with uniform bounds
 depending on $r,k, e_0,e_1$.
Next
$$
|(D_t^2 \phi)(x)|\le 2C \int \int_0^{2z} |\phi (y)|dy \mu_t(dz)
$$
 for all $x$, which by Cauchy-Schwartz inequality does not exceed
$$
2C  \|\phi \|_{L_{2,f_k}} \int \sqrt{\int_0^{2z} f_k^2(y)\, dy}
\mu_t (dz)
 \le Cc(k) \|\phi \|_{L_{2,f_k}}(1+(E^{k+1/2}, \mu_0)).
 $$
 Hence writing the solution to the Cauchy problem for equation \ref{eq5.13}
 with $\phi_r \in L_{2,f_k}$ as a perturbation series \eqref{eqA16}
 with respect to perturbation $D_t^2$ one represents the solution as
 the sum $\phi_1^t+\phi_2^t$ with
 $$
 \|\phi_1^t\|_{L_{2,f_k}} \le c(k,e_0,e_1,r)\|\phi_r\|_{L_{2,f_k}}
 $$
 and
 $$
 \|\phi_2^t\|_{C(\R_+)} \le (1+(E^{k+1/2}, \mu_0))\kappa(C,k,e_0,e_1,r)\|\phi_r\|_{L_{2,f_k}}
 $$
 so that
 $$
 \|\phi_2^t\|_{L_{2,f_k}} \le (1+(E^{k+1/2}, \mu_0))\kappa(C,k,e_0,e_1,r)\|\phi_r\|_{L_{2,f_k}}
 $$
 whenever $k> 1/2$. In particular for these $k$
 $$
 \|U^{t,r} \|_{L^{1,0}_{2,f_k}} \le (1+(E^{k+1/2}, \mu _0))
 \kappa (C,k,e_0,e_1,r).
 $$
As $(\xi_s(\mu_0;x;.),g)=(\delta_x,U^{0,s}g)$, this implies
\eqref{eq5.20} by \eqref{eq1.81}. The evolution $U^{t,r}$ in the
space $L^{2,0}_{2,f_k}$ is analyzed quite similarly.

We conclude with the following analog of Proposition \ref{prop5.5},
whose proof follows from Proposition \ref{prop5.7} by the same
argument as Proposition \ref{prop5.5} follows from Proposition
\ref{prop5.4}.

\begin{proposition}\label{prop5.8}
Under the condition (C3) for $k>1/2$ and $m=1,2$
\begin{equation}\label{eq5.22}
\sup_{s \leq t} \| \mu_s(\mu^1_0) -
\mu_s(\mu^2_0)\|_{(L^{m,0}_{2,f_k})'}
 \leq \kappa (C,t,k,e_0,e_1)(1+E^{1+k/2},\mu^1_0+\mu_0^2)
  \|\mu^1_0-\mu^2_0\|_{(L^{m,0}_{2,f_k})'}
\end{equation}
\end{proposition}

\section{The rate of convergence in the LLN} \label{sec6}

{\it Proof of Theorem \ref{th1}.} Recall that $\mu_t (Y)$ means the
solution to equation \eqref{eq1.5} with initial data $\mu_0=Y$ given
by Proposition \ref{prop2.1} with a $\beta \geq 2$. We shall write
shortly $Y_t=\mu_t(Y)$ so that $T_tF(Y)=F(\mu_t(Y))=F(Y_t)$.
 For a function $F(Y) = (g, Y^{\otimes m})$
with
 $g \in C^{sym}_{(1 +
E)^\otimes,\infty}(X^m)$, $m \geq 1$, and $Y = h \delta_\x$ one has
\begin{equation}\label{eq6.1}
T_tF(Y) - T^h_tF(Y) =\int_0^t T^h_{t-s} (L_h - \LC)T_s F(Y)\, ds.
\end{equation}
As $T_tF(Y)=(g,Y_t^{\otimes m})$, Propositions \ref{prop4.2} and
\ref{prop4.3} yield
$$
\delta T_tF(Y;x) = m \int_{X^m}g(y_1, y_2, \ldots, y_m) \xi_t(Y ;
x;dy_1) Y_t^{\otimes(m-1)}(dy_2 \cdots dy_m),
$$
 and
\begin{align}\label{eq6.2}
& \delta^2 T_tF(Y;x,w)
= m\int_{X^m} g(y_1, y_2, \ldots, y_m) \eta_t(Y ; x,w;dy_1) Y_t^{\otimes(m-1)}(dy_2 \cdots dy_m)\nonumber  \\
& + m(m-1)\int_{X^m} g(y_1,y_2, \ldots, y_m) \xi_t (Y ; x;dy_1)
\xi_t (Y ; w;dy_2) Y_t^{\otimes(m-2)}(dy_3 \cdots dy_m).
\end{align}
Let us estimate the difference $(L_h - \LC)T_tF(Y)$ using
\eqref{eq3.5} (with $T_tF$ instead of $F$). Let us analyze only the
more weird second term in \eqref{eq3.5}, as the first one is
analyzed similar, but much simpler. We are going to estimate
separately the contribution to the last term of \eqref{eq3.5} of the
first and second term in \eqref{eq6.2}.

Assume that the condition (C1) or (C2) holds and a $k \geq 1$ is chosen.
Note that the norm and the first moment
$(E,\cdot)$ of $Y + sh(\delta_y - \delta_{x_i} - \delta_{x_j})$
 do not exceed respectively the norm and the first moment of $Y$. Moreover, for $s\in [0,1]$,
 $h>0$ and $x_i,x_j,y \in X$ with $E(y)=E(x_i)+E(x_j)$ one has
 \begin{align*}
(E^k,Y &+ sh(\delta_y - \delta_{x_i} - \delta_{x_j})) = (E^k,Y) + sh(E(x_i) + E(x_j))^k - hE^k(x_i) - hE^k(x_j) \\
 & \leq (E^k,Y) + hc(k)(E^{k-1}(x_i) E(x_j) + E(x_i) E^{k-1}(x_j))
 \end{align*}
 with a constant $c(k)$ depending only on $k$. Consequently by Proposition \ref{prop5.2}
 \begin{align*}
 & \| \eta_t (Y + sh(\delta_y - \delta_{x_i} - \delta_{x_j});x,w;\cdot)\|_{1 + E^k}
 \leq \kappa (C,k,t,e_0,e_1)(1+E(w)) \\
 & \{1+E^{k+1}(x)+[(E^{k+1},Y)+hc(k)(E^k(x_i)E(x_j)+E(x_i)E^k(x_j))](1+E^2(x)) \\
  & + [(E^{k+3},Y)+hc(k)(E^{k+2}(x_i)E(x_j)+E(x_i)E^{k+2}(x_j)](1+E(x))\}+...,
 \end{align*}
 where by dots is denoted the similar term with $x$ and $w$ interchanging their places.
 Hence the contribution to the last term of \eqref{eq3.5} of the first term in \eqref{eq6.2} does not exceed
 \begin{align*}
 & \kappa (C,t,k,m,e_0, e_1) \|g\|_{(1 + E^k)^{\otimes m}}(1+E^k,Y)^{m-1}h^3 \sum_{i \neq j}
 (1+E(x_i)+E(x_j))^2 \{1+E^{k+1}(x_i)+E^{k+1}(x_j) \\
 & +[(E^{k+1},Y)+hc(k)(E^k(x_i)E(x_j)+E(x_i)E^k(x_j))](1+E^2(x_i)+E^2(x_j)) \\
  & + [(E^{k+3},Y)+hc(k)(E^{k+2}(x_i)E(x_j)+E(x_i)E^{k+2}(x_j)](1+E(x_i)+E(x_j))\}.
 \end{align*}
 Dividing this sum into two parts, where $E(x_i) \geq E(x_j)$ and respectively vice versa,
 and noting that by the symmetry it is enough to estimate only the first part,
  allows to estimate the contribution to the last term of \eqref{eq3.5} of the first term
  from \eqref{eq6.2} by
 \begin{align*}
  & \kappa (C,t,k,m,e_0, e_1) \|g\|_{(1 + E^k)^{\otimes m}}(1+E^k,Y)^{m-1}h^3 \sum_{i \neq j} \\
 & \{1+E^{k+3}(x_i)+(1+E^4(x_i))
 [(E^{k+1},Y)+hc(k)E^k(x_i)E(x_j)] \\
  & + (1+E^3(x_i))[(E^{k+3},Y)+hc(k)E^{k+2}(x_i)E(x_j)]\}.
 \end{align*}
The main term in this expression (obtained by ignoring the terms with $hc(k)$) is estimated by
 $$
 \kappa \|g\|_{(1 + E^k)^{\otimes m}} ( 1 + E^k,Y)^{m-1}h
  [(1+E^{k+3},Y)+(1+E^4,Y)(E^{k+1},Y)+(1+E^3,Y)(E^{k+3},Y)],
  $$
  where the first two terms in the square bracket can be estimated by the last one, because
  $$
  (E^4,Y)(E^{k+1},Y)\le 2(E^2,Y)(E^{k+3},Y).
  $$
  It remains to observe that
the terms with $hc(k)$ are actually subject to the same bound, as for instance
$$
 h^4 \sum_{i \neq j}E^k(x_i)E(x_j)(1+E^4(x_i)) \leq h^2(E^k+E^{k+4},Y)(E,Y)
  \leq c(k) h(E^{k+3},Y)(E,Y)^2.
 $$
Consequently, the the contribution to the last term of \eqref{eq3.5}
of the first term in \eqref{eq6.2} does not exceed
\begin{equation}
\label{eq6.3}
 h\kappa(C,t,k,m,e_0,e_1) \|g\|_{(1+E^k)^{\otimes m}}
(1+E^k,Y)^{m-1} (1+E^{k+3},Y)(1+E^3,Y).
\end{equation}
Turning to the contribution of the second term from \eqref{eq6.2}
observe that again by Proposition \ref{prop5.2}
\begin{align*}
 & \| \xi_t (Y + sh(\delta_y - \delta_{x_i} - \delta_{x_j});x;\cdot)\|_{1 + E^k}
 \leq \kappa (C,k,t,e_0,e_1) \\
 & \{1+E^k(x)+(1+E(x)[(E^{k+1},Y)+hc(k)(E^k(x_i)E(x_j)+E(x_i)E^k(x_j))]\},
 \end{align*}
 so that the contribution of the second term from \eqref{eq6.2} does
 not exceed
 \begin{align*}
 & \kappa (C,t,k,m,e_0, e_1) \|g\|_{(1 + E^k)^{\otimes m}}(1+E^k,Y)^{m-2}h^3 \sum_{i \neq j}
 (1+E(x_i)+E(x_j)) \{1+E^k(x_i)+E^k(x_j) \\
 &
 +(1+E(x_i)+E(x_j))[(E^{k+1},Y)+hc(k)(E^k(x_i)E(x_j)+E(x_i)E^k(x_j))]\}^2,
 \end{align*}
 which again by
 dividing this sum into two parts, where $E(x_i) \geq E(x_j)$ and respectively vice versa,
 reduces to
\begin{align*}
 & \kappa (C,t,k,m,e_0, e_1) \|g\|_{(1 + E^k)^{\otimes m}}(1+E^k,Y)^{m-2}h^3 \sum_{i \neq j}
 (1+E(x_i)) \{1+E^k(x_i) \\
 &
 +(1+E(x_i)+E(x_j))[(E^{k+1},Y)+hc(k)E^k(x_i)E(x_j)]\}^2.
 \end{align*}
 This is again estimated by \eqref{eq6.3}. It
 follows now from
\eqref{eq6.1} and Proposition \ref{prop3.1} that
\[
\|T_tF - T^h_t F\|_{C_{(1 + E^{k+3},\cdot)(1+E^3,\cdot)(1+E^k,\cdot)^{m-1}}(\MC_{h \delta}^{e_0e_1}(X))}
  \leq h \kappa (C,t,k,m,e_0,e_1) \|g\|_{(1 + E^k)^{\otimes m}},
  \]
which is the same as \eqref{eq2.5}. The proof of \eqref{eq2.6} is
quite the same. It only uses Proposition \ref{prop5.4} instead of
Proposition \ref{prop5.2}.

\section{Auxiliary Estimates}\label{sec7}

The main technical ingredient in the proof of a weak form of CLT
(convergence for fixed times, stated in Theorems \ref{th2}-
\ref{th4}) is given by the following corollary to Theorem \ref{th1}.
\begin{proposition}\label{prop7.1}
Under the assumptions of Theorem \ref{th1} let $g_2$ be a symmetric
continuous function on $X^2$. Then for any $k\ge 1$
\begin{equation}\label{eq7.1}
\sup_{s \leq t} \quad \left|\E \left( g_2, \left( \frac{Z^h_s(Z^h_0)
- \mu_s(\mu_0)}{\sqrt{h}} \right)^{\otimes 2} \right)\right| =
\sup_{s \leq t}|\E (g_2, (F_s^h (Z_0^h,\mu_0))^{\otimes 2})|
=\sup_{s \leq t}|(U^{h;0,s}_{fl} (g_2,.))(F_0^h)|
\end{equation}
 does not exceed the expression
$$
 \kappa (C,t,k,e_0,e_1) \|g_2\|_{(1 +
E^k)^{\otimes 2}(X^2)}(1+(E^{k+3},Z^h_0+\mu_0))^2
  \left( 1 + \left\| \frac{Z^h_0 -\mu_0}{\sqrt{h}} \right\|^2_{1 + E^k} \right)
$$
for any $k\ge 1$ under the condition (C1) or (C2) and the expression
$$
 \kappa (C,t,k,e_0,e_1)
\|g_2\|_{C^{2,sym}_{(1 + E^k)^{\otimes
2}(X^2)}}(1+(E^{k+4},Z^h_0+\mu_0))^3
  \left( 1 + \left\| \frac{Z^h_0 -\mu_0}{\sqrt{h}} \right\|^2_{\MC^1_{1 + E^k}} \right)
$$
 for any $k\ge 0$ under the condition (C3) with a constant $\kappa (C,t,k,e_0,e_1)$.
\end{proposition}

{\it Proof.}
One has
\[
\E \left( g_2, \left( \frac{Z^h_t(Z^h_0) - \mu_t(\mu_0)}{\sqrt{h}}
\right)^{\otimes 2} \right) =\E \left( g_2, \left(
\frac{Z^h_t(Z^h_0) - \mu_t(Z^h_0)}{\sqrt{h}} \right)^{\otimes 2}
\right) + \left( g_2, \left( \frac{\mu_t(Z^h_0) -
\mu_t(\mu_0)}{\sqrt{h}} \right)^{\otimes 2} \right)
\]
\begin{equation}
\label{eq7.2}
 + 2 \E \left(
g_2, \frac{Z^h_t(Z^h_0) - \mu_t(Z^h_0)}{\sqrt{h}}
  \otimes  \frac{\mu_t(Z^h_0) - \mu_t(\mu_0)}{\sqrt{h}} \right)
\end{equation}
The first term can be rewritten as
\begin{align*}
  \frac{1}{h} \E & \bigl(g_2, (Z^h_t (Z^h_0))^{\otimes 2} - (\mu_t(Z^h_0))^{\otimes 2} \\
  & + \mu_t(Z^h_0) \otimes (\mu_t(Z^h_0) - Z^h_t (Z^h_0)) + (\mu_t(Z^h_0) - Z^h_t(Z^h_0)) \otimes \mu_t (Z^h_0)\bigr).
  \end{align*}
  Under the condition (C1) or (C2) this term can be estimated by
  $$
  \kappa (C,r,e_0, e_1) \|g_2\|_{(1+E^k)^{\otimes 2}}(1+(E^{k+3},Z^h_0))
 (1+(E^k,Z^h_0))(1+(E^3,Z^h_0))
 $$
 $$
  \leq \kappa (C,r,e_0, e_1) \|g_2\|_{(1+E^k)^{\otimes 2}} (1 + (E^{k+3}, Z^h_0))^2,
  $$
  due to Theorem \ref{th1} and \eqref{eq5.2}. The second term is estimated by
  $$
  \|g_2\|_{(1+E^k)^{\otimes 2}}(1+(E^{k+1},\mu_0+Z^h_0))\left\| \frac{Z^h_0 -\mu_0}{\sqrt{h}} \right\|^2_{1 + E^k}
  $$
  by \eqref{eq2.3}, and the third term by the obvious combination
   of these two estimates completing the proof for cases (C1) and (C2). The case (C3) is considered
   analogously.
   Namely, the first term in the representation \eqref{eq7.2} is
   again estimated by Theorem \ref{th1}, and to estimate the second
   term one uses \eqref{eq5.15} instead of \eqref{eq2.3} and the observation that
   $$
   |(g_2,\nu^{\otimes 2})| \le \sup_{x_1} \left| (1+E^k(x_1))^{-1}
   \int {\partial g_2 \over \partial x_1}(x_1,x_2) \nu(dx_2)\right| \Vert \nu
   \Vert_{\MC^1_{1+E^k}}
 $$
 $$
 \le \sup_{x_1,x_2} \left| (1+E^k(x_1))^{-1}(1+E^k(x_2))^{-1}
   {\partial^2 g_2 \over \partial x_1 \partial x_2}(x_1,x_2) \right| \Vert \nu
   \Vert_{\MC^1_{1+E^k}}^2
   \le \Vert g_2 \Vert_{(C^{2,sym}_{(1+E^k)^{\otimes 2}}} \Vert \nu
   \Vert_{\MC^1_{1+E^k}}^2.
   $$

   Though the estimates of Proposition \ref{prop7.1} are sufficient
   to prove Theorem \ref{th2}, in order to prove the semigroup convergence
   from Theorem \ref{th4} one needs a slightly more general
   estimate, which in turn requires a more general form of LLN,
   than presented in Theorem \ref{th1}. We shall give now these two
   extensions.

 \begin{proposition}\label{prop7.2}
 The estimates on the r.h.s. of \eqref{eq2.5} and \eqref{eq2.6}
 remain valid, if one the l.h.s. on takes a more general expression,
 namely
 $$
 \sup_{s\le t} |T_t^h(GFH)(Y)-G(Y_t)F(Y_t)T_t^hH(Y)|,
 $$
 where $F(Y)$ is as in Theorem \ref{th1} and both $G$ and $H$ are
 cylindrical functionals of the form \eqref{eq2.12} with $f \in C^2(\R^d)$ and
 all $\phi_j$, $j=1,...,n$, belonging to $C_{1+E^k}(X)$ and $C^{2,0}_{1+E^k}(X)$
 respectively in cases (C1)-(C2) and (C3) (with a
 constant $C$ depending on the corresponding norms of $\phi_j$).
 \end{proposition}

 {\it Proof.} As
 $$
 T_t^h(GFH)(Y)-G(Y_t)F(Y_t)T_t^hH(Y)
 $$
 $$
 =T_t^h(GFH)(Y)-(GFH)(Y_t)
 + (GF)(Y_t)(H(Y_t)-T_t^hH(Y)),
 $$
 it is enough to consider the case without a function $H$ involved.
 And in this case looking through the proof of Theorem \ref{th1}
 above one sees that it generalizes straightforwardly to give the
 result required.

 \begin{proposition}\label{prop7.3}
  The estimates of Proposition \ref{prop7.1} remain valid if instead
  of \eqref{eq7.1} one takes a more general expression
  \begin{equation}\label{eq7.3}
 \sup_{s \leq t} \left|\E [(g_2, F_s^h (Z_0^h,\mu_0))G(F_s^h(Z_0^h,\mu_0))]\right|
   =\sup_{s \leq t} |[U^{h;0,s}_{fl} ((g_2,.)G)](F_0^h)|,
\end{equation}
where $G$ is as in the previous Proposition (with a constant $C$
again depending on the norms of
 $\phi_j$ in the representation of $G$ as a cylindrical function of the form \eqref{eq2.12}).
\end{proposition}

{\it Remark.} Let us stress for clarity that $U^{h;0,s}_{fl}
((g_2,.)G)$ in \eqref{eq7.3} means the result of the evolution
$U^{h;0,s}_{fl}$ applied to the function of $Y$ given by
$(g_2,Y^{\otimes 2})G(Y)$.

 {\it Proof.} It is again obtained by a straightforward
generalization of the proof of Proposition \ref{prop7.1} given above
using Proposition \ref{prop7.2} instead of Theorem \ref{th1}.

The main technical ingredient in the proof of the functional CLT
(stated as Theorems \ref{th5}- \ref{th6}) is given by the following

\begin{proposition}\label{prop7.4}
Under condition (C3) for any $k>1/2$
 \begin{equation}
 \label{eq7.4}
 \sup_{s \leq t}\E \| F_s^h (Z_0^h,\mu_0) \|^2_{(L^{2,0}_{2,f_k})'}
   \leq \kappa (C,t,k,e_0,e_1) (1+(E^{k+5},Z^h_0+\mu_0))^3
  (1+ \| F_0^h \|^2_{(L^{2,0}_{2,f_k})'}).
\end{equation}
\end{proposition}

{\it Proof.} The idea is to represent the l.h.s. of \eqref{eq7.4} in
the form of the l.h.s. of \eqref{eq7.1} with an appropriate function
$g_2$. Using the notation $\tilde \nu(x) =\int_x^{\infty} \nu (dy)$
from the introduction for a finite (signed) measure $\nu$ on $\R_+$
(and setting $\tilde \nu (x)=0$ for $x<0$) one has
$$
\FC (\tilde \nu)
 =\frac{1}{\sqrt {2\pi}} \int_0^{\infty} e^{-ipy} (\int_y^{\infty}\nu (dx))dy
 =\frac{1}{\sqrt {2\pi}} \int_0^{\infty}\nu (dx) \int_0^x e^{-ipy}\, dy,
 $$
 so that for $f_k(x)=1+x^k$
$$
\FC (f_k\tilde \nu)
 =(1+(i\frac{\partial}{\partial p})^k)\FC (\tilde \nu)
 =\frac{1}{\sqrt {2\pi}} \int_0^{\infty}\nu (dx) \int_0^x (1+y^k)e^{-ipy}\, dy.
 $$
 Applying \eqref{eq1.82} yields
 $$
 \|\nu\|^2_{(L^{2,0}_{2,f_k}(\R_+))'}= \|f_k\tilde \nu\|^2_{H^{-1}(\R)}
  =\frac{1}{2\pi}\int_{-\infty}^{\infty}\left| \int_0^{\infty} \nu (dx) \int_0^x
  (1+y^k)e^{-ipy} \, dy\right|^2 \frac{dp}{1+p^2}
  $$
  \begin{equation}\label{eq7.5}
  =\frac{1}{2\pi}\int_0^{\infty} \int_0^{\infty} \theta_k(x,y) \nu (dx)\nu (dy)
  \end{equation}
  with
  \begin{equation}\label{eq7.6}
  \theta_k(x,y)=\int_{-\infty}^{\infty}
  \left(\int_0^x (1+z^k)e^{-ipz}\, dz
  \int_0^y (1+w^k)e^{ipw}\, dw \right) \frac{dp}{1+p^2}.
  \end{equation}
  Clearly
  $$
  \theta_k(x,y)|_{x=0}=\theta_k(x,y)|_{y=0}=0.
  $$
  Moreover
  $$
  \frac{\partial \theta_k}{\partial x}(x,y)
   =\int_{-\infty}^{\infty}[(1+x^k)e^{-ipx}\int_0^y
   (1+w^k)e^{ipw}dw] \frac{dp}{1+p^2}
   $$
   so that
   $$
 |\frac{\partial \theta_k}{\partial x}(x,y)|\le
 c(k)(1+x^k)(1+y^{k+1}),
 $$
 and
  $$
  |\frac{\partial^2 \theta_k}{\partial x \partial y}(x,y)|
   =\left|\int_{-\infty}^{\infty}(1+x^k)e^{-ipx}
   (1+y^k)e^{ipy} \frac{dp}{1+p^2}\right|\leq c(k)(1+x^k)(1+y^k).
   $$
   Since
$$
  \frac{\partial^2 \theta_k}{\partial x^2}(x,y)
   =\int_{-\infty}^{\infty}[kx^{k-1}e^{-ipx}\int_0^y
   (1+w^k)e^{ipw}dw] \frac{dp}{1+p^2}
   $$
   $$
   -\int_{-\infty}^{\infty}[(1+x^k)e^{-ipx}\int_0^y
   (1+w^k)(ip)e^{ipw}dw] \frac{dp}{1+p^2},
   $$
   and using integration by parts in the second term yields also
   $$
 |\frac{\partial^2 \theta_k}{\partial x^2}(x,y)|
  \leq c(k)[(1+x^k)(1+y^k)+(1+x^{k-1})(1+y^{k+1})].
  $$
  Consequently
  $\theta_k +\bar \theta_k\in C^{2,sym}_{(1+E^{k+1})^{\otimes 2}}$.
  Therefore, using \eqref{eq7.5} for $\nu=F_s^h(Z_0^h,\mu_0)$
  implies that in order to estimate the l.h.s. of \eqref{eq7.4} one needs to
  estimate the l.h.s. of \eqref{eq7.1} with $g_2=\theta_k$ given by
  \eqref{eq7.6}.

  Though a direct
  application of Proposition \ref{prop7.1} does not give the result
  we need, only a slight modification is required. Namely, representing
  \eqref{eq7.1} in form \eqref{eq7.2} we estimate the first
  term precisely like in the proof of Proposition \ref{prop7.1} and the second
  term that now equals
  $$
  \left\| \frac{\mu_t(Z_0^h)-\mu_t(\mu_0)}{\sqrt h} \right\|^2_{(L^{2,0}_{2,f_k})'}
  $$
  can be estimated using Proposition \ref{prop5.8} by
  $$
  \kappa (c,t,k,e_0,e_1)(1+E^{k+1},\mu_0+Z^h_0)
  \| \frac{Z_0^h-\mu_0}{\sqrt h} \|^2_{(L^{2,0}_{2,f_k})'}.
  $$
  Estimating the third term in \eqref{eq7.2} again by the
  combination of the estimates of the first two terms yields
  \eqref{eq7.4}.



\section{CLT: Proof of Theorems \ref{th2} - \ref{th6}}\label{sec8}

  {\it Proof of Theorem \ref{th2}.} Recall that we denoted by $U^{h;t,r}_{fl}$ the backward propagator
   corresponding to the process $F^h_t=(Z^h_t - \mu_t) / \sqrt{h}$.
   By \eqref{eq3.10}, the l.h.s. of \eqref{eq2.10} can be written as
  $$
  \sup_{s \leq t} \left| (U^{h;0,s}_{fl}(g,.))(F^h_0)- (U^{0,s}g, F^h_0)\right|
 = \sup_{s \leq t} \int_0^s [U^{h;0,\tau}_{fl}  (\Lambda^h_{\tau} - \Lambda_{\tau})
  U^{\tau,s}(g,.)] \, d\tau \, (F^h_0).
  $$
  As by \eqref{eq3.9}
  \begin{align*}
  (\Lambda^h_{\tau}&  - \Lambda_{\tau})(U^{{\tau},s}g,\cdot)(Y)
  = \frac{\sqrt{h}}{2} \int \int \int (U^{{\tau},s} g(y) - U^{{\tau},s} g(z_1) - U^{{\tau},s}g(z_2))
   K(z_1, z_2;dy)Y(dz_1)Y(dz_2) \\
  & - \frac{\sqrt{h}}{2} \int \int  (U^{{\tau},s}g(y) - 2U^{{\tau},s} g(z)) K(z,z;dy)(\mu_t+\sqrt{h}Y)(dz)
  \end{align*}
  (note that the terms with the second and third variational
  derivatives in  \eqref{eq3.9} vanish here, as we apply it to a
  linear function),
  the required estimate follows from Proposition \ref{prop7.1}.

{\it Proof of Theorem \ref{th3}}.
  Substituting the function $\Phi_{f_t}^{\phi_1^t,...,\phi_n^t}$ of form \eqref{eq2.12} (with two times continuously
differentiable $f_t$) with a given initial condition
$\Phi_r(Y)=\Phi_{f_r}^{\phi_1^r,...,\phi_n^r}(Y)$ at $t=r$ in the
equation $\dot F_t=-\Lambda_t F_t$ yields
$$
\frac{\partial f_t}{\partial t} + \frac{\partial f_t}{\partial
x_1}(\dot \phi_1^t,Y)+... +\frac{\partial f_t}{\partial x_n}(\dot
\phi^t_n,Y)
$$
$$
=-\frac{1}{2} \int \int \int \sum_{j=1}^n \frac{\partial
f_t}{\partial x_j}(\phi_j^t,
\delta_y-\delta_{z_1}-\delta_{z_2})K(z_1,z_2;dy)(Y(dz_1)\mu_t(dz_2)+\mu_t(dz_1)Y(dz_2))
$$
$$
-\frac{1}{4} \int \int \int \sum_{j,l=1}^n \frac{\partial^2
f_t}{\partial x_j \partial x_l} (\phi_j^t\otimes \phi_l^t,
(\delta_y-\delta_{z_1}-\delta_{z_2})^{\otimes
2})K(z_1,z_2;dy)\mu_t(dz_1)\mu_t(dz_2)
$$
with $f_t(x_1,...,x_n)$ and all its derivatives evaluated at the
points $x_j=(\phi^t_j,Y)$ (here and in what follows we denote by dot
the derivative $d/dt$ with respect to time). This equation is
clearly satisfied whenever
\begin{equation}
\label{eq8.1}
\dot f_t(x_1,...,x_n) =-\sum_{j,k=1}^n\Pi (t, \phi_j^t,\phi_k^t)
\frac{\partial ^2 f_t}{\partial x_j \partial x_k}(x_1,...,x_n)
\end{equation}
and
$$
\dot \phi_j^t(z) =-\int \int (\phi_j^t(y)-\phi_j^t(z)-\phi_j^t(w))
 K(z,w;dy) \mu_t(dw)=-\Lambda_t \phi_j^t(z)
$$
with $\Pi$ given by \eqref{eq2.16}. Consequently
\begin{equation}
\label{eq8.2}
OU^{t,r} \Phi_r(Y)=\Phi_t(Y)=(\UC^{t,r}
f_r)((U^{t,r}\phi_1^r,Y),...,(U^{t,r}\phi_n^r,Y)),
\end{equation}
where $\UC^{t,r}f_r=\UC ^{t,r}_{\Pi}f_r$ is defined as the resolving
operator to the (inverse time) Cauchy problem of equation
\eqref{eq8.1} (it is well defined as \eqref{eq8.1} is just a
spatially invariant second order evolution), the resolving operator
$U^{t,r}$ is constructed in Sections 4,5, and
 $$
 \Pi(t,\phi_j^t,\phi_k^t)=\Pi (t, U^{t,r}\phi_j^r, U^{t,r}\phi_k^r).
 $$
  All statements of
Theorem \ref{th3} follows from the explicit formula \eqref{eq8.2},
the semigroup property of the solution to finite-dimensional
equation \eqref{eq8.1} and Propositions \ref{prop5.4},
\ref{prop5.7}.

{\it  Proof of Theorem \ref{th4}.}
 The first statement is obtained by a straightforward modification of our proof of Theorem
 \ref{th2} above, where one has to use Proposition \ref{prop7.3}
 instead of its particular case Proposition \ref{prop7.1} and to
 note that all terms in
 formula \eqref{eq3.9} (that unlike the linear case now become relevant) depend
  at most quadratically on $Y$, because for a function $\Phi$ of form \eqref{eq2.12}
\[
\delta \Phi (Y;x) = \sum_{j=1}^n \frac{\partial f}{\partial x_j}
\phi_j(x),
\]
\begin{equation}
\label{eq8.21}
 \delta ^2 \Phi (Y;x,y) = \sum_{i,j=1}^n
\frac{\partial^2 f}{\partial x_j \partial x_i} \phi_j(x) \phi_i(y),
\end{equation}
where the derivatives of $f$ are evaluated at the points
$x_j=(\phi^t_j,Y)$.

    The second statement follows by
  the usual approximation of a general $\Phi$ by those given by
  \eqref{eq2.12} with smooth $f$.

{\it  Proof of Theorem \ref{th5}.}  The characteristic function of
$\Phi_t^h$ is
$$
g_{t_1,...,t_n}^h(p_1,...,p_n)
 ={\bf E} \exp \{ \sum_{j=1}^n (\phi_j, F^h_{t_j} (Z^h_0,\mu_0))\}
 =U_{fl}^{h;0,t_1}\Phi_1... U_{fl}^{h;t_{n-2},t_{n-1}}\Phi_{n-1}
 U_{fl}^{h;t_{n-1},t_n}\Phi_n(F_0^h),
 $$
 where
 $\Phi_j(Y)= \exp \{ip_j (\phi_j,Y)\}$. Let us show that it converges
  to the characteristic function
\begin{equation}\label{eq8.3}
g_{t_1,...,t_n}(p_1,...,p_n)
 =OU^{0,t_1}\Phi_1... OU^{t_{n-2},t_{n-1}}\Phi_{n-1}
 OU^{t_{n-1},t_n}\Phi_n(F_0)
 \end{equation}
 of a Gaussian random variable. For $n=1$ it follows from Theorem
 \ref{th4}. For $n>1$ one can write
 $$
 g_{t_1,...,t_n}^h(p_1,...,p_n)-g_{t_1,...,t_n}(p_1,...,p_n)
 $$
 \begin{equation}\label{eq8.4}
  =\sum_{j=1}^n U_{fl}^{h;0,t_1}\Phi_1...U_{fl}^{h;t_{j-2},t_{j-1}}
  \Phi_{j-1} (U_{fl}^{h;t_{j-1},t_j}-OU^{t_{j-1},t_j}) \Phi_j
  OU^{t_j,t_{j+1}}...OU^{t_{n-1}t_n} \Phi_n.
  \end{equation}
  By Theorem \ref{th4} we know that for any $j=1,...,n$
  $$
  \Psi_j^{p_j,...,p_n}(Y)=(U_{fl}^{h;t_{j-1},t_j}-OU^{t_{j-1},t_j}) \Phi_j
  OU^{t_j,t_{j+1}}...OU^{t_{n-1}t_n} \Phi_n(Y)
  $$
  converge to zero as $\sqrt h$ as $h\to 0$ uniformly on $Y$ from
  bounded domains of $\MC^1_{1+E^k}$.
  We have to show that
  \begin{equation}
  \label{eq8.5}
  U_{fl}^{h;t_{j-2},t_{j-1}} \Phi_{j-1}\Psi_j(Y)
  ={\bf E}^h_Y(\Phi_{j-1} \Psi_j (Y_{t_{j-1}}))
  \end{equation}
  tends to zero, where ${\bf E}_Y^h$ is of course the expectation
  with respect to the fluctuation process started in $Y$ at time
  $t_{j-2}$. The last expression can be written as
  \begin{equation}
  \label{eq8.6}
  {\bf E}^h_Y(({\bf 1}_{\{\|Y_{t_{j-1}}\|_{(L^{2,0}_{2,f_{k+2}})'}\le K\}}
 \Phi_{j-1} \Psi_j )(Y_{t_{j-1}}))
  +{\bf E}^h_Y(({\bf 1}_{\{\|Y_{t_{j-1}}\|_{(L^{2,0}_{2,f_{k+2}})'}> K\}}
 \Phi_{j-1} \Psi_j )(Y_{t_{j-1}})).
 \end{equation}
  For $Y$ from a bounded subset of $(L^{2,0}_{2,f_{k+2}})'$ the second
  term can be made arbitrary small by choosing large enough $K$ due
  to Proposition \ref{prop7.4}. Due to the natural continuous
  inclusion $C_{f_k}^{m,0} \subset L^{m,0}_{2,f_{k+\alpha}}$,
  $m=1,2$ , $\alpha >1/2$ one gets by duality a continuous
  projection $(L_{2,f_k}^{2,0})' \mapsto \MC^2_{f_{k-\alpha}}
  \subset \MC^1_{f_{k-\alpha}}$ for $k>1/2$, $\alpha \in (1/2,k)$.
  Hence a bounded set in $(L_{2,f_{k+2}}^{2,0})'$ is also bounded in
  $\MC^1_{f_{k+2-\alpha}}$, so that there
  $\Phi_{j-1}\Psi_j(Y_{t_{j-1}})$ is small of order $\sqrt h$,
  implying
  that the first term in \eqref{eq8.6} is small. Consequently
  expression \eqref{eq8.5} tends to zero uniformly for $Y$ from
  bounded domain of $(L_{2,f_{k+2}}^{2,0})'$, $k>1/2$. This implies that
  all terms in \eqref{eq8.4} tend to zero as $h\to 0$.

  It remains to check that \eqref{eq8.3} is given by \eqref{eq2.15},
  which is done by induction in $n$
  using Theorem
 \ref{th3} and an obvious explicit formula
 $$
\UC^{t,r}f(x)=\exp \{i \sum_{j=1}^n p_j x_j -\sum_{j,k=1}^n p_jp_k
\int_t^r \Pi (s, \phi_j^s, \phi_k^s) \, ds \}
$$
for the solution of the
 Cauchy problem of the diffusion equation \eqref{eq8.1}
with $f(x)=\exp \{i\sum_{j=1}^n p_jx_j \}$. For instance,
 $$
 (OU^{t_{n-1},t_n}\Phi_n)(Y)
 =(\UC^{t_{n-1},t_n}f_n)(U^{t_{n-1},t_n}\phi_n,Y)
 $$
 $$
 =\exp \{ip_n(U^{t_{n-1},t_n}\phi_n,Y)
  -p_n^2\int_{t_{n-1}}^{t_n} \Pi (s,
  U^{s,t_n}\phi_n,U^{s,t_n}\phi_n) \, ds \}
 $$
 where $f_n(x)=\exp \{ip_nx\}$, and hence
 $$
 OU^{t_{n-2},t_{n-1}}(\Phi_{n-2} OU^{t_{n-1},t_n}\Phi_n)(Y)
 $$
 $$
 =\exp
 \{i(p_{n-1}U^{t_{n-2},t_{n-1}}\phi_{n-1}+p_nU^{t_{n-2},t_n}\phi_n,Y)
 -p_n^2\int_{t_{n-2}}^{t_n} \Pi (s, U^{s,t_n}\phi_n,U^{s,t_n}\phi_n) \, ds \}
 $$
 $$
 \times \exp \{-\int_{t_{n-2}}^{t_{n-1}}
 [p_{n-1}^2 \Pi (s, U^{s,t_{n-1}}\phi_{n-1}, U^{s,t_{n-1}}\phi_{n-1})
 +2p_{n-1}p_n \Pi (s, U^{s,t_{n-1}}\phi_{n-1}, U^{s, t_n}\phi_n) ]\, ds \}.
 $$
The proof is complete.

{\it  Proof of Theorem \ref{th6}.}

(i) Notice first that applying Dynkin's formula to the Markov
process $Z^h_t$ one finds that for a $\phi \in C_{1+E^k}(X)$
\[
M^h_{\phi}(t)=(\phi, Z_t^h)-(\phi,Z^h_0)-\int_0^t (L_h
(\phi,.))(Z^h_s) ds
\]
is a martingale, since all three terms here are integrable, due to
formula \eqref{eq3.61} and the assumption $Z_0^h \in
\MC_{1+E^{k+5}}$. Hence $(\phi,F^h_t)$ is a semimartingale and
\[
(\phi,F^h_t)=\frac{M^h_{\phi}}{\sqrt h}+V_{\phi}^h(t)
\]
with
\[
V_{\phi}^h(t)= \frac{1}{\sqrt h}
 [\phi, Z_0^h)+\int_0^t (L_h (\phi,.))(Z^h_s)
ds -(\phi,\mu_t)]
\]
is the canonical representation of the semimartingale $(\phi,F^h_t)$
into the sum of a martingale and a predictable process of bounded
variation that is also continuous and integrable. (It implies, in
particular that $(\phi,F^h_t)$ belongs to the class of special
semimartingales.)

As we know already the convergence of finite dimensional
distributions, to prove (i) one has to show that the distribution on
the Skorohod space of c\`adl\`ag functions of the semimartingale
$(\phi, F_t^h)$ is tight, which according to Aldous-Rebolledo
Criterion (see e.g. \cite{Et}, \cite{Re}, we cite the formulation
from \cite{Et}) amounts to showing that given a sequence of $h_n \to
0$ as $n\to \infty$ and a sequence of stopping times $\tau_n$
bounded by a constant $T$ and an arbitrary $\epsilon>0$ there exist
$\delta>0$ and $n_0>0$ such that
$$
\sup_{n\ge n_0} \sup_{\theta \in [0,\delta]} P\left[ |V^{(n)}
(\tau_n+\theta)-V^{(n)}(\tau_n)|>\epsilon\right] \le \epsilon,
$$
and
$$
\sup_{n\ge n_0} \sup_{\theta \in [0,\delta]} P\left[ |Q^{(n)}
(\tau_n+\theta)-Q^{(n)}(\tau_n)|>\epsilon\right] \le \epsilon,
$$
where $V^n(t)$ is a shorter notation for $V^{h_n}_{\phi}$ and
$Q^n(t)$ is the quadratic variation of the martingale
$M^{h_n}_{\phi}(t)$. Notice that it is enough to show the tightness
of $(\phi,F_t^h)$ for a dense subspace of the test functions $\phi$.
Thus we can and will consider now only the bounded $\phi$.

 To get a required
estimate for $V^n(t)$ observe that by \eqref{eq3.61}
\[
\frac{d}{dt} V^n(t)=\frac{1}{\sqrt h}
 [(L_h (\phi,.))(Z^h_s)
 -(\phi, \dot {\mu}_t)]
 =- \frac{ \sqrt h}{2}\int \int [\phi (y)-2\phi (z)]K(z,z;dy)Z^h_t(dz)
 \]
 \[
 + \frac{1}{2 \sqrt h} \int \int \int [\phi (y) -
 \phi (z_1) - \phi (z_2)]
 K(z_1, z_2;dy) [Z^h_t(dz_1)Z^h_t(dz_2)-\mu_t(dz_1)\mu_t(dz_2)].
 \]
 The first term here is clearly uniformly bounded for $h\to 0$, and
 the second term can be written as
\begin{equation}
  \label{eq8.61}
\frac{1}{2} \int \int \int [\phi (y) -
 \phi (z_1) - \phi (z_2)]
 K(z_1, z_2;dy) [F^h_t(dz_1)Z^h_t(dz_2)+\mu_t(dz_1)F^h_t(dz_2)].
 \end{equation}
 Applying Doob's maximal inequality to the martingale
\[
(\phi, F_t^h)-(\phi,F^h_0)-\int_0^t (\Lambda^h_s (\phi,.))(Z^h_s) ds
\]
in combination with Proposition \ref{prop7.4} shows that
  \eqref{eq8.61} can be
made bounded with an arbitrary
 small probability, implying the required estimate for $V^n(t)$.

 Let us estimate the
quadratic variation by the same arguments as in \cite{Ko3}. Namely,
as the process $(\phi, F_t^{h})$ for each $h$ is the sum of a
differentiable process and a pure jump process, both having locally
finite variation, its quadratic variation coincides with that of
$M^{h_n}_{\phi}(t)$ and is known to equal the sum of the squares of
the sizes of all its jumps (see e.g. Theorem 26.6 in \cite{Kal}), so
that
$$
Q^{(n)}(t)-Q^{(n)}(\tau)=\sum_{s\in [\tau,t]} (\phi,
F^{h_n}_s-F^{h_n}_{s-})^2
 ={1 \over h}\sum_{s\in [\tau,t]} (\phi,
Z^{h_n}_s-Z^{h_n}_{s-})^2.
$$
As each jump of $Z^h_s$ is the change of $h\delta_x+h\delta_y$ to
$h\delta_{x+y}$ one concludes that
$$
|Q^{(n)}(t)-Q^{(n)}(\tau)| \le h \sup |\phi||N_t-N{\tau}|
$$
with $N_t$ denoting the number of jumps on the interval $[0,t]$. By
the L\'evy formula for Markov chains (see e.g. \cite{Br}) the
process $N_t-\int_0^t a(Z_s^h) \, ds$ is a martingale, where $a(Y)$
denotes the intensity of jumps at $Y$, given by \eqref{eq2.4}.
Hence, using the optional sampling theorem and \eqref{eq2.4} implies
that
\[
\E(N_t-N_{\tau})=\E\int_{\tau}^t a(Z_s^h) \, ds
 \le 3Ch^{-1} e_0(e_1+e_0)\E (t-\tau),
\]
and consequently
$$
\E|Q^n(t)-Q^n(\tau)| \le 3C \|\phi\|e_0(e_1+e_0)\theta
$$
uniformly for all $<t-\theta<\tau<t$. Hence by Chebyshev inequality
the required estimate for $Q^n$ follows.

(ii) By Theorem \ref{th5} the limiting process is uniquely defined
whenever it exists. Hence one only needs to prove the tightness of
the family of normalized fluctuations $F^h_t$. Again due to the
existence of finite dimensional limits and general convergence
theorems (see either a result of \cite{Mi} specially designed to
show convergence in Hilbert spaces, or a more general result on
convergence of a complete separable metric space valued processes in
\cite{EK} or \cite{Et}), to prove tightness it is enough to
establish the following compact containment condition: for every
$\epsilon>0$ and $T>0$ there exists $K>0$ such that for any $h$
$$
P(\sup_{t\in [0,T]}\|F^h_t\|_{(L^{2,0}_{2,f_k})'}>K) \le \epsilon.
$$
To this end, let us introduce a regularized square root function
$R$, i.e. $R(x)$ is an infinitely smooth increasing function
$\R_+\mapsto \R_+$ such that $R(x)=\sqrt x$ for $x>1$, and the
corresponding "regularized norm" functional on $(L^{2,0}_{2,f_k})'$:
$$
G(Y)=R((Y,Y)_{(L^{2,0}_{2,f_k})'})=R((\theta_k,Y\otimes Y)),
$$
where $\theta_k$ is given by \eqref{eq7.5} (see Proposition
\ref{prop7.4}). By Dynkin's formula one can conclude that the
process
$$
M_t=G(F^h_t)-G(F_0^h)-\int_0^t \Lambda_s^h G(F_s^h) \, ds
$$
is a martingale whenever all terms in this expression have finite
expectations. (Note that we use here a more general than usual
version of Dynkin's formula with a time dependent generator; the
reduction of time nonhomogeneous case to the standard situation by
including time as an additional coordinate of a Markov process under
consideration is explained e.g. in \cite{Fr}.) Expectation of
$G(F^h_t)$ is bounded by Proposition \ref{prop7.4}. Moreover, taking
into account \eqref{eq8.21} and the fact that
$R^{(k)}(s)=O(s^{(1/2)-k})$ for $s\ge 1$, one sees from formulas
\eqref{eq3.9} and \eqref{eq2.8} that $\Lambda_s^h G(F_s^h)$ grows at
most quadratically in $F_s^h$, which again by Proposition 7.4
implies the uniform boundedness of the expectation of this term.
Applying to $M_t$ Doob's maximal inequality yields the required
compact containment completing the proof of the theorem.

\section{Three lemmas}\label{sec9}

We present here three general (not connected to each other) analytic
facts used in the main body of the paper. Recall that classes
$C^1(\MC_f(X))$ were defined in the introduction.

\begin{lemma}\label{lem9.1}
(i) If $F\in C^1(\MC_f(X))$ and $Y,\xi \in \MC_f(X)$, then
\begin{equation}
\label{eq1.9}
 F(Y + \xi) - F(Y) = \int_0^1 (\delta F(Y + s \xi ;
\cdot), \xi)\, ds.
\end{equation}

(ii) If $F \in C^2(\MC_{f, \phi}(X))$ or $F \in C^3(\MC_{f,
\phi}(X))$, the following Taylor expansion holds respectively:
\begin{align}
\label{eq1.10}
 (a) \hspace{1cm} & F(Y + \xi) - F(Y) = (\delta F(Y;
\cdot), \xi) + \int_0^1 (1 - s) (\delta^2 F(Y + s \xi ; \cdot,
\cdot), \xi \otimes \xi)\, ds,
  \nonumber
 \\
(b) \hspace{1cm} &F(Y + \xi) - F(Y)
 = (\delta F(Y; \cdot), \xi) + \frac{1}{2} (\delta^2 F(Y; \cdot, \cdot), \xi \otimes \xi) \nonumber \\
 & + \frac{1}{2} \int_0^1 (1-s)^2 (\delta^3 F(Y + s \xi; \cdot, \cdot, \cdot), \xi^{\otimes 3})\, ds.
 \end{align}

(iii) If $t \mapsto \mu_t \in \MC_f(X)$ is continuous in the
$*$-weak topology of $\MC_f(X)$ and is continuously differentiable
in the
 $*$-weak topology of $\MC_{\phi}(X)$, then for any $F \in C^1(\MC_{f, \phi}(X))$, $\phi \leq $f,
\begin{equation}
\label{eq1.11}
\frac{d}{dt} F(\mu_t) = (\delta F(\mu_t; \cdot), \dot{\mu}_t).
 \end{equation}
\end{lemma}

\noindent {\it Proof.} (i)
 Using the representation
 $$
F(Y + s(\delta_x+\delta_y)) - F(Y)
 =F(Y + s\delta_x) - F(Y)
  +\int_0^s \delta F(Y+s \delta_x+h\delta_y;y)\, dh
 $$
 for arbitrary points $x,y$
 and the uniform continuity of $\delta F(Y+s \delta_x+h\delta_y;y)$
 in $s,h$ allows to deduce from \eqref{eq1.83} the existence of the limit
$$
\lim_{s \rightarrow 0_+} \frac{1}{s} (F(Y + s(\delta_x+\delta_y)) -
F(Y))=\delta F(Y;x)+\delta F(Y;y).
$$
Extending similarly to the arbitrary number of points one obtains
\eqref{eq1.9} for $\xi$ being an arbitrary finite sum of the Dirac
measures $\delta_{x_1}+...+\delta_{x_n}$.

 Assume now that $\xi \in \MC_f(X)$ and $\xi_k \to \xi$ as $k\to
\infty$ $\star$-weakly in $\MC_f(X)$, where all $\xi_k$ are finite
sums of Dirac measures. We are going to pass to the limit $k\to
\infty$ in the equation \eqref{eq1.9} written for $\xi_k$. As $F\in
C(\MC_f)$ one has
$$
F(Y+\xi_k)-F(Y) \to F(Y+\xi)-F(Y), \quad k\to \infty.
$$
Next, the difference
$$
\int_0^1 (\delta F(Y + s \xi_k ; \cdot), \xi_k)\, ds
 -\int_0^1 (\delta F(Y + s \xi ;
\cdot), \xi)\, ds
$$
can be written as
$$
\int_0^1 (\delta F(Y + s \xi_k ; \cdot), \xi_k-\xi)\, ds
 +\int_0^1 (\delta F(Y + s \xi_k ;
\cdot)- \delta F(Y + s \xi ; \cdot), \xi)\, ds.
$$
The second term tends to zero, because by our assumption the
variational derivation $\delta F$ maps $\MC_f(X)$ continuously to
$C_{f,\infty}(X)$. The first term tends to zero, because $\xi_k \to
\xi$ weakly and the family of functions $\delta F(Y+s\xi_k;.)$ is
compact in $C_{f,\infty}(X)$ (which is again due to the assumed
continuity of the derivation $\delta F$).

 Statement (ii) is
straightforward from the usual Taylor expansion. Turning to (iii)
observe that
$$
\frac{d}{dt} F(\mu_t)=\lim_{h\to 0} \frac {1}{h}
(F(\mu_{t+h})-F(\mu_t)),
$$
which by (i) and the assumed continuous differentiability can be
written as
$$
\lim_{h\to 0} \int_0^1 \left(\delta F (\mu_t +s
(\mu_{t+h}-\mu_t);.), \frac {1}{h} \int_0^h \dot \mu_{t+\tau}
d\tau\right) \, ds.
$$
We want to show that it equals the r.h.s. of \eqref{eq1.11}. We have
\[
 \int_0^1 \left(\delta F (\mu_t +s
(\mu_{t+h}-\mu_t);.), \frac {1}{h} \int_0^h \dot \mu_{t+\tau}
d\tau\right) \, ds - (\delta F(\mu_t; \cdot), \dot{\mu}_t)
\]
\[
=\left(\delta F (\mu_t;.), \frac {1}{h} \int_0^h \dot \mu_{t+\tau}
d\tau - \dot{\mu}_t\right)
 +\int_0^1 \left(\delta F (\mu_t +s
(\mu_{t+h}-\mu_t);.)-\delta F(\mu_t; \cdot), \frac {1}{h} \int_0^h
\dot \mu_{t+\tau} d\tau\right) \, ds.
\]
The first term here tends to zero as $h\to 0$ by the weak continuity
of $\dot \mu_t$. The second term tends to zero, because the family
of measures
 \[
\mu_{t+h}-\mu_t=\frac{1}{h}\int_0^h \dot \mu_{t+\tau} d\tau
\]
is bounded and hence compact in the $\star$-weak topology of
$\MC_{\phi}(X)$.

\begin{lemma}\label{lem9.2}
Suppose $S$ is a compact subset of a linear topological space $Y$
(we are interested in the case when $Y$ is a topological dual to a
Banach space equipped with its $\star$-weak topology) and $Z_t$ is a
Markov process on $S$ specified by its Feller semigroup $\Psi_t$ on
$C(S)$ with a bounded generator $A$. Let  $\Omega_t(z) = (z -
\xi_t)/a$ be a family of linear transformation on $Y$, where $a$ is
a positive constant and $\xi_t$, $t\geq 0$, is a differentiable
curve in $Y$. Let
\[
\Omega_{[0,T]}(S)=\cup_{t\in [0,T]} \Omega_t (S)
\]
for a $T>0$.
 Then $Y_t = \Omega_t(Z_t)$, $t\in [0,T]$, is a Markov process in $\Omega_{[0,T]}(S)$ for any
 $T>0$ with the dynamics of averages (propagator)
 \[
 U^{s,t}f(y)=\E_{s,y} f(Y_t)
 \]
 given by the formula
\begin{equation}
\label{eq9.10}
 U^{s,t}f(y)=\Omega^{-1}_s \Psi_{t-s}\Omega_t f(y)
\end{equation}
for any $f\in C(\Omega_{[0,T]}(S))$, $t\le T$,
 where $\Omega_tf(y) = f (\Omega_t(y))$. Moreover, if such a
 function $f$ is uniformly continuously differentiable in the
 direction $\dot \xi_t$, i.e. if the limit
 \[
 \lim_{\tau \to 0} \frac{1}{\tau}\left(f(\Omega_{t+\tau}(y))
  -f(\Omega_t(y))\right)=-\frac{1}{a}(\nabla_{\dot \xi_t}f)(\Omega_t(y))
  =\nabla_{\dot \xi_t}f(\Omega_t(y))
  \]
  exists and is uniform in $\Omega_{[0,T]}(S)$, then for all
   $s\le t$
\begin{equation}
\label{eq9.11}
  \frac{d}{dt}U^{s,t}f=U^{s,t}\Lambda_tf,
\end{equation}
where the operator $\Lambda_t$ is given by the formula
\begin{equation}\label{eq9.13}
\Lambda_t f = \Omega^{-1}_t A \Omega_t f - \frac{1}{a} \nabla_{\dot
\xi_t}f.
\end{equation}
\end{lemma}

{\em Proof.} Formula \eqref{eq9.10} follows from the definitions of
$\Psi_t$ and $\Omega_t$. Formulas \eqref{eq9.11}, \eqref{eq9.12}
follow by differentiating \eqref{eq9.10} using the product rule and
taking into account that the derivative $\nabla f$ is supposed to be
uniform.

{\it Remark.} Similarly, using the identity
 \[
\Omega_t^{-1} \nabla_{\dot \xi_t}\Omega_t =a^{-1}\nabla_{\dot
\xi_t},
 \]
one shows that
\begin{equation}
\label{eq9.12}
  \frac{d}{ds}U^{s,t}f=-\Lambda_s U^{s,t} f,
\end{equation}
holds for $s=t$. However, to extend this to $s<t$ one needs some
additional assumptions on the smoothness of the semigroup $\Psi_t$.

\begin{lemma}\label{lem9.3} \cite{Ko3}
Let $Y$ be a measurable space and
 the mapping $t\mapsto \mu_t$ from $[0,T]$ to ${\cal M}(Y)$ is
 continuously differentiable in the sense of the norm in ${\cal M}(Y)$
 with a (continuous) derivative $\dot \mu_t=\nu_t$. Let $\sigma_t$ denote a
 density of $\mu_t$
 with respect to its total variation $|\mu_t|$, i.e. the
 class of measurable functions (equivalence is defined as the a.s. equality with respect to
 the measure $|\mu_t|$) taking three values $-1,0,1$ and such that
 $\mu_t=\sigma_t|\mu_t|$ and $|\mu_t|=\sigma_t\mu_t$
 almost surely with respect to $|\mu_t|$. Then there exists a
 measurable function $f_t(x)$ on $[0,T]\times Y$ such that
 $f_t$ is a representative
 of class $\sigma_t$ for any $t \in [0,T]$ and
 $$
 \Vert \mu_t \Vert=\Vert \mu_0 \Vert +\int_0^t ds \int_Yf_s(y)\nu_s(dy).
$$
\end{lemma}

We refer for a proof to the Appendix of \cite{Ko3} noting only that
$f_t$ could be chosen as such a representative of $\sigma_t$, which
is at the same time a representative of the class of the densities
of $\nu_t^s$ with respect to its total variation measure
 $|\nu_t^s|$, where $\nu_t^s$ is a singular part of $\nu_t$ in its Lebesgue decomposition
 with respect to $|\mu_t|$.

 \appendix

 \section{On the evolutions with integral generators}\label{appendix}
 Here we present an analytic study of evolutions with integral
 generators that are obtained as certain perturbations of positivity
 preserving evolutions. As always, it is assumed that $X$ is a
 locally compact space (though this assumption is used only in Theorem \ref{thA2},
 other statement being valid for arbitrary topological spaces).

 We shall start with the problem
 \begin{equation}\label{eqA1}
 \dot{u}_t(x) = A_t u_t(x) = \int u_t(z) \nu_t(x;dz) - a_t(x) u_t(x), \quad u_r(x) = \phi(x), \quad t \geq r \geq 0,
 \end{equation}
 where $\phi$ and $a_t$ are given measurable functions on $X$ such that
 $a_t$ is non-negative and locally bounded in $t$ for each
 $x$, $\nu_t(x,\cdot)$ is a given family of finite (non-negative) measures on
 $X$ depending measurably on $t \geq 0$, $x \in X$, and such that
 $\sup_{t \in [0,T]} \| \nu_t(x,\cdot)\|$ is bounded for arbitrary $T$ and $x$.

 Clearly equation \eqref{eqA1} is formally equivalent to the integral equation
 \begin{equation}\label{eqA2}
 u_t(x) = I^r_\phi(u)_t = e^{-(\xi_t(x) - \xi_r(x))} \phi(x) + \int_r^t e^{-(\xi_t(x) - \xi_s(x))} L_s u_s(x)ds,
 \end{equation}
 where $\xi_t(x) = \int_0^t a_s(x)ds$ and $L_tv(x) = \int v(z) \nu_t(x,dz)$.

 We shall look for the solutions of \eqref{eqA2} in the class of
 functions $u_t(x)$, $t \geq r$, that are continuous in $t$ (for each $x$),
 measurable in $x$ and such that the integral in the expression for $Lu_s$
 is well defined in the Lebesgue sense. Basic obvious observation about
 \eqref{eqA2} is the following: the iterations of the mapping $I^r_\phi$
 form \eqref{eqA2} are connected with the partial sums
 $$
 S^{t,r}_m \phi = \left[ e^{-(\xi_t-\xi_r)}
 + \sum_{l=1}^m \int_{r \leq s_l \leq \cdots \leq s_1 \leq t} e^{-(\xi_t - \xi_{s_1})}
 L_{s_1} \cdots L_{s_{l-1}}e^{- (\xi_{s_{l-1}} - \xi_{s_l})} L_{s_l}e^{-(\xi_{s_l} - \xi_r)}
  ds_1 \cdots ds_l\right] \phi
 $$
(where $e^{-\xi_t}$ designates the operator of multiplication by
$e^{-\xi_t(x)}$) of the perturbation series solution $S^{t,r} =
\lim_{m \rightarrow \infty} S^{t,r}_m$ to \eqref{eqA2} by
\begin{equation}\label{eqA3}
(I^r_\phi)^m(\phi)_t = S^{t,r}_{m-1} \phi + \int_{r \leq s_m \leq
\cdots \leq s_1 \leq t} e^{-(\xi_t - \xi_{s_1})}L_{s_1} \cdots
L_{s_{m-1}}e^{-(\xi_{s_{m-1}} - \xi_{s_m})}L_{s_m}\phi ds_1 \cdots
ds_m.
\end{equation}
\begin{lemma}\label{lemA1}
Suppose
\begin{equation}\label{eqA4}
A_t \psi(x) \leq c \psi(x), \quad t \in [0,T],
\end{equation}
for a strictly positive measurable function $\psi$ on $X$ and a constant $c = c(T)$. Then
\begin{equation}\label{eqA5}
(I^r_\psi)^m(\psi)_t \leq \left( 1 + c(t-r) + \cdots + \frac{1}{m!} c^m(t-r)^m\right) \psi
\end{equation}
for all $0 \leq r \leq t \leq T$, and consequently $S^{t,r} \psi$ is well defined as
 a convergent series for each $t,x$ and
\begin{equation}\label{eqA6}
S^{t,r} \psi(x) \leq e^{c(t-r)} \psi(x).
\end{equation}
\end{lemma}

{\it Proof.} This is given by induction in $m$. Suppose \eqref{eqA5} holds for $m$. Since \eqref{eqA4} implies
$$L_t \psi(x) \leq (c + a_t(x)) \psi(x) = (c + \dot{\xi}_t(x)) \psi(x),$$
it follows that
\begin{align*}
(I^r_\psi)^{m+1}&(\psi)_t \leq e^{-(\xi_t(x) - \xi_r(x))} \psi(x) \\
& + \int_r^t e^{-(\xi_t(x) - \xi_s(x))}(c + \dot{\xi}_s(x))
 \left(1 + c(s-r) + \cdots + \frac{1}{m!} c^m(s-r)^m\right) \psi(x)\, ds.
\end{align*}
Consequently, as
$$
\int_r^t e^{-(\xi_t - \xi_s)} \dot{\xi}_s \frac{1}{l!} (s-r)^l \, ds = \frac{1}{l!} (t-r)^l - \frac{1}{(l-1)!} \int_r^t
e^{-(\xi_t - \xi_s)}(s-r)^{l-1}\, ds
$$
 for $l > 0$, it remains to show that
\begin{align*}
\sum_{l=1}^m& c^l \left[ \frac{1}{l!} (t-r)^l - \frac{1}{(l-1)!} \int_r^t e^{-(\xi_t - \xi_s)}(s-r)^{l-1}ds\right] + \sum_{l=0}^m c^{l+1} \frac{1}{l!} \int_r^t e^{-(\xi_t - \xi_s)}(s-r)^lds \\
& \leq c(t-r) + \cdots + \frac{1}{(m+1)!} c^{m+1}(t-r)^{m+1}.
\end{align*}
But this holds, because the l.h.s. of this inequality equals
$$\sum_{l=1}^m \frac{c^l}{l!} (t-r)^l + \frac{c^{m+1}}{m!} \int_r^t e^{-(\xi_t - \xi_s)}(s-r)^m ds.$$

The following corollary plays an important role in the analysis of
Section 5.

\begin{lemma}\label{lemA2}
Suppose $A_t \psi \leq c \psi + \phi$ for positive functions $\phi$ and $\psi$ and all $t \in [0,T]$. Then
$$
S^{t,r} \psi \leq e^{c(t-r)}[\psi + \int_r^tS^{\tau,t}\, d\tau
\phi].$$
\end{lemma}

{\it Proof.} Using \eqref{eqA5} yields
\begin{align*}
I^r_\psi(\psi)_t & \leq (1 + c(t-r))\psi + \int_r^t e^{-(\xi_t - \xi_s)} \phi \, ds, \\
( I^r_\psi)^2(\psi)_t & \leq \left(1 + c(t-r) + \frac{c^2}{2}(t-r)^2\right)\psi
 + \int_r^t e^{-(\xi_t - \xi_s)}(1+c(s-r)) \phi \, ds \\
& + \int_r^t e^{-(\xi_t - \xi_s)} L_s \int_r^s e^{-(\xi_s -
\xi_\tau)} \phi \, d \tau \, ds,
\end{align*}
etc, and hence
$$
( I^r_\psi)^m (\psi)_t \leq e^{c(t-r)}
 \left[ \psi + \int_r^t e^{-(\xi_t - \xi_s)} \phi \, ds
  + \int_r^t e^{-(\xi_t - \xi_s)}
L_s \int_r^s e^{-(\xi_s - \xi_\tau)} \phi \, d \tau ds + \cdots
\right]
 $$
  $$
  =e^{c(t-r)}
 \left[ \psi + \int_r^t d\tau \left( e^{-(\xi_t - \xi_{\tau})}
  + \int_{\tau}^t e^{-(\xi_t - \xi_s)}
L_s e^{-(\xi_s - \xi_{\tau})} \, ds + \cdots \right)\phi \right]
$$
and the proof is completed by noting that
$$
S^{t,r} \psi = \lim_{m \rightarrow \infty} S^{t,r}_{m-1} \psi
 \leq \lim_{m \rightarrow \infty} (I^r_\psi)^m (\psi)_t.
 $$
 The existence of the solutions to \eqref{eqA1} and \eqref{eqA2} can be easily established now.

 \begin{proposition}\label{propA1}
 Under the assumptions of Lemma \ref{lemA1} the following holds.

 (i) For an arbitrary $\phi \in B_\psi$ the perturbation series
 $S^{t,r} \phi = \lim_{m \rightarrow \infty} S^{t,r}_m \phi$ is absolutely
 convergent for all $t,x$, the function $S^{t,r} \phi$ solves
 \eqref{eqA2} and represents its minimal solution
 (i.e. $S^t \phi \leq u$ point-wise for any other solution $u$ to \eqref{eqA2}),
  and $S^{t,r} \phi (x)$ tends to
  $S^{\tau,r}\phi(x)$ as $t \rightarrow \tau$ uniformly on any set where both $a_t$ and $\psi$ are bounded.

 (ii) The family $S^{t,r}$ form a propagator in $B_\psi(X)$ with the norm
\begin{equation}\label{eqA7}
\| S^{t,r}\|_\psi \leq c^{C(t-r)}.
\end{equation}
\end{proposition}

{\it Proof.}  Applying Lemma \ref{lemA1} separately to the positive
and negative part of $\phi$ one obtains the convergence of series
$S^{t,r} \phi$ and the estimate \eqref{eqA7}. Clearly $S^{t,r} \phi$
satisfies \eqref{eqA2} and it is minimal, as any solution $u$
of this equation satisfies the equation $u_t = (I^r_\phi)^m(u)_t$ and
hence (due to \eqref{eqA3}) also the inequality $u_t \geq S^{t,r}_{m-1} \phi$.

The continuity of $S^{t,r}$ in $t$ follows from the formula
$$
S^{t,r}\phi-S^{\tau,r}\phi=(e^{-(\xi_t-\xi_{\tau})}-1)e^{-(\xi_{\tau}-\xi_r)}\phi
$$
\begin{equation} \label{eqA8}
+ \int_r^{\tau} (e^{-(\xi_t - \xi_\tau)}-1)e^{-(\xi_\tau - \xi_s)}
L_s S^{s,r} \phi \,ds
 + \int_\tau^t e^{-(\xi_t - \xi_s)} L_s
S^{s,r}\phi \, ds
\end{equation}
for $r \leq \tau \leq t$.

At last, once the convergence of the series $S^{t,r}$ is proved,
the propagator (or Chapman-Kolmogorov) equation \eqref{eq1.12} follows
from simple standard manipulations with integrals that we omit.

For the application to time non-homogeneous stochastic processes
 one needs actually equation \eqref{eqA1} in inverse time, i.e. the problem
\begin{equation}\label{eqA9}
\dot{u}_t(x) = - A_t u_t(x) = - \int u_t(z) \nu_t(x;dz) + a_t(x)u_t(x), \quad u_r(x) = \phi(x), \quad 0 \leq t \leq r,
\end{equation}
with the corresponding integral equation taking the form
\begin{equation}\label{eqA10}
u_t(x) = I^r_\phi(u)_t = e^{\xi_t(x) - \xi_r(x)} \phi(x) + \int_t^r e^{\xi_t(x) - \xi_s(x)}
L_s u_s(x)ds.
\end{equation}

All the statements of Proposition \ref{propA1} (and their proofs)
obviously hold for the perturbation series $S^{t,r}$ constructed
from \eqref{eqA10}, with the same estimate \eqref{eqA7}, but with
the backward propagator equation \eqref{eq1.12} holding for $t \leq
s \leq r$ with $S$ instead of $U$.

To get a strong continuity of $S^{t,r}$ one usually needs a second
bound for $A_t$. In particular, the following holds.

\begin{proposition}\label{propA2}
Suppose now that two measurable functions $\psi_1, \psi_2$ on $X$ are given both
 satisfying \eqref{eqA4} and such that (i) $0 < \psi_1 < \psi_2$,
  (ii) $a_t$ is bounded on any set where $\psi_2$ is bounded,
  (iii) $\psi_1 \in B_{\psi_2, \infty}(X)$.
Then $S^{t,r}$, $t\leq r$ (constructed above for \eqref{eqA9}, \eqref{eqA10})
is a strongly continuous family of operators in
$B_{\psi_2, \infty}(X)$.
\end{proposition}

{\it Proof.} By Proposition \ref{propA1} $S^{t,r}$ are bounded in $B_{\psi_2}(X)$.
 Moreover, as $S^{t,r}\phi$ tends to $\phi$
uniformly on the sets where $\psi_2$ is bounded, it follows that
$$
\| S^{t,r} \phi - \phi\|_{\psi_2} \rightarrow 0
$$
for any $\phi \in B_{\psi_1}(X)$, and hence also for any
 $\phi \in B_{\psi_2,\infty}(X)$, since $B_{\psi_1}(X)$ is dense in
$B_{\psi_2,\infty}(X)$.

\begin{theorem}\label{thA1}
Under the assumptions of Proposition \ref{propA2} assume additionally
 that $\psi_1, \psi_2$ are continuous, $a_t$ is a continuous mapping
  $t \mapsto C_{\psi_2/\psi_1,\infty}$ and $L_t$ is a continuous mapping
   from $t$ to bounded operators $C_{\psi_1} \mapsto C_{\psi_2,\infty}$.
    Then $B_{\psi_1}$ is an invariant core for the propagator $S^{t,r}$ in the sense that
\begin{align}\label{eqA11}
A_r \phi & = \lim_{t \rightarrow r, t\leq r} \frac{S^{t,r} \phi - \phi}{r-t}
 = \lim_{s\rightarrow r, s \geq r} \frac{S^{r,s} \phi - \phi}{s-r}, \nonumber \\
\frac{d}{ds} S^{t,s} \phi & = S^{t,s} A_s \phi,
 \quad \frac{d}{ds} S^{s,r} \phi = - A_s S^{s,r} \phi, \quad t < s < r,
\end{align}
for all $\phi \in B_{\psi_1}(X)$, with all these limit existing in the
 Banach topology of $B_{\psi_2,\infty}(X)$. Moreover, $C_{\psi_1}$
  and $C_{\psi_2,\infty}$ are invariant under $S^{t,r}$, so that
   $C_{\psi_1}$ is an invariant core of the strongly continuous propagator
    $S^{t,r}$ in $C_{\psi_2,\infty}$. In particular, if $a_t, L_t$ do not
     depend on $t$, then $A$ generates a strongly continuous semigroup
  on $C_{\psi_2,\infty}$ with $C_{\psi_1}$ being an invariant core.
\end{theorem}

{\it Proof.} The differentiability of $S^{t,r} \phi(x)$ for each $x$ follows from
\eqref{eqA8} (better to say its time reversal version).
Differentiating equation \eqref{eqA10} one sees directly that $S^{t,r} \phi$
satisfies \eqref{eqA9} and al required formulas hold point-wise.
To show that they hold in the topology of $B_{\psi_2,\infty}$ one needs
to show that the operators $A_t(\phi)$ are continuous as functions from
$t$ to $B_{\psi_2,\infty}$ for each $\phi \in B_{\psi_1}$.
But this follows directly from our continuity assumptions on $a_t$ and $L_t$.

To show that the space $C_{\psi_1}$ is invariant (and this wold obviously imply all
 other remaining statements), we shall approximate $S^{t,r}$ by the evolutions
  with bounded intensities. Let $\chi_x$ be a measurable function $X \mapsto [0,1]$
   such that $\chi_n(x) = 1$ for $\psi_2(x) \leq n$ and $\chi_n(x) = 0$ for $\psi_2(x) \geq n+1$.
    Denote $\nu^n_t(x,dz) = \chi_n(x) \nu_t(x,dz)$, $a^n_t = \chi_n a_t$,
     and let $S^{t,r}_n$ (respectively $A^n_t$) denote the propagators constructed as
  in Proposition \ref{propA2} (respectively  the operators from \eqref{eqA1}) but with
   $\nu^n_t$ and $a^n_t$ instead of $\nu_t$ and $a_t$. Then the propagators $S^{t,r}_n$
    converge strongly in the Banach space $B_{\psi_2,\infty}$ to the propagator $S^{t,r}$.
     One can deduce this fact from a general statement on the convergence of propagators
   (see e.g. \cite{Ma1}), but a direct proof is even simpler. Namely, as $S^{t,r}$ and
    $S^{t,r}_n$ are uniformly bounded, it is enough to show the convergence for the elements
  $\phi$ of the invariant core $B_{\psi_1}$. For such a $\phi$ one has
\begin{equation}\label{eqA12}
(S^{t,r} - S^{t,r}_n)(\phi) = \int_t^r \frac{d}{ds} S^{t,s} S^{s,r}_n \phi \, ds
 = \int_t^r S^{t,s} (A_s - A^n_s) S^{s,r}_n \phi \, ds,
\end{equation}
where \eqref{eqA11} was used. As by invariance $S^{s,r}_n \phi \in B_{\psi_1}$, it follows
 that $(A_s - A^n_s) S^{s,r}_n \phi \in B_{\psi_2}$ and tends to zero in the form of
  $B_{\psi_2}$, as $n \rightarrow \infty$, and hence the r.h.s. of \eqref{eqA12} tends
   to zero in $B_{\psi_2}$, as $n \rightarrow \infty$.

To complete the proof it remains to observe that as the generators of $S^{t,r}_n$ are bounded,
 the corresponding semigroups preserves continuity (as they can be constructed as the convergent
  exponential series). Hence $S^{t,r}$ preserves the continuity as well, as $S^{t,r} \phi$
   is a (uniform) limit of continuous functions.

{\it Remark.} Choosing $a_t = \| \nu_t(x,\cdot)\|$ and $\psi_1 = 1$ above yield a
pure analytic construction of a strongly continuous propagator for a non-homogeneous
jump type process. A more familiar probabilistic approach can be found e.g.
in \cite{Ch} (at least for the homogeneous case).

For our purposes we need a perturbed equation \eqref{eqA9}, namely the equation
\begin{equation}\label{eqA13}
\dot{u}_t = -(A_t - B_t)u_t, \quad u_r = \phi, \quad 0 \leq t \leq r,
\end{equation}
where $B_t$ are bounded operators in $C_{\psi_1}$, and its dual equation on measures, whose weak form is
\begin{equation}\label{eqA14}
\frac{d}{dt} (g, \xi_t) = ((A_t - B_t)g, \xi_t) \quad \xi_0 = \xi, \quad 0 \leq t \leq r,
\end{equation}
i.e. has to hold for some class of test functions $g$. Motivated by the standard observation that formally
 equation \eqref{eqA13} is equivalent to the integral equation

\begin{equation}\label{eqA15}
u_t = S^{t,r} \phi - \int_t^\tau S^{t,s} B_s u_sds,
\end{equation}
whose solution $u_t = U^{t,r} \phi$ one expects to obtain through the perturbation series
\begin{equation}\label{eqA16}
U^{t,r} \phi = S^{t,r} \phi - \int_t^r S^{t,s} B_s S^{s,r} ds
 + \int_{t \leq s_1 \leq s_2 \leq r} S^{t, s_1} B_{s_1} S^{s_1, s_2} B_{s_2} S^{s_2,r} ds_1 ds_2 + \cdots,
\end{equation}
one arrives at the following result.

\begin{theorem}\label{thA2}
Under the assumptions of Theorem \ref{thA1} suppose that $\psi_2(x) \rightarrow \infty$
as $x \rightarrow \infty$ and that a strongly continuous family of bounded operators
 $B_t :C_{\psi_2} \mapsto C_{\psi_1}$ is given. Then

(i) series \eqref{eqA16} is absolutely convergent in $C_{\psi_2}(X)$
for any $\phi \in C_{\psi_2}(X)$ so that
$$
\Vert U^{t,r}\Vert_{C_{\psi_2}(X)} \leq \Vert
S^{t,r}\Vert_{C_{\psi_2}(X)} \exp \{(r-t) \sup_{t\le s\le r} \Vert
B_s \Vert_{C_{\psi_2}(X)} \},
$$
 and defines a strongly
continuous backward propagator $U^{t,r}$ in $C_{\psi_2,\infty}(X)$ with $C_{\psi_1}$ being its
 invariant core (so that the analogues of \eqref{eqA11} hold);

(ii) the operator $V^{r,s} = (U^{s,r})^*$ form a weakly continuous propagator in $\MC_{\psi_2}$ yielding a unique (weakly
continuous) solution to the Cauchy problem \eqref{eqA14} in the sense that it holds for all $g \in C_{\psi_1}(X)$;

(iii) if $f$ is an arbitrary continuous function tending to zero as $x \rightarrow \infty$,
then the operators $V^{r,s} = (U^{r,s})^*$ are strongly continuous in the norm of
$\MC_{\psi_2 f}$ and solves a strong version of \eqref{eqA14} with derivative
taken in the norm topology of $\MC_{\psi_1 f}$.

(iv) at last, if a family $A_t^\omega, B_t^\omega$ of operators are given satisfying
 all the above conditions for each $\omega$ from an interval and such that
  $A^\omega_t - B^\omega_t$ depend strongly continuous on $\omega$ as operators
  $C_{\psi_1} \mapsto C_{\psi_2,\infty}$, then the corresponding resolving
  operators $U^{s,r}$ in $C_{\psi_2,\infty}$ depend strongly continuous
  on $\omega$ and their adjoint operators $V^{r,s}$ depend weakly continuous on $\omega$ in $\MC_{\psi_2}$.
\end{theorem}

{\it Proof.} (i) \eqref{eqA16} converges, because $B_t$ are bounded.
Other statements then follow directly from the corresponding facts about $S^{t,r}$.

(ii) The operators $V^{r,s}$ are weakly continuous in $\MC_{\psi_2}(X)$ just because they
 are adjoint to strongly continuous
operators in $C_{\psi_2,\infty}$. Next, the analogue of the third equation
 in \eqref{eqA11} for $U^{t,r}$ is the equation
$$
\frac{d}{dr} U^{s,r} g = U^{s,r} (A_r - B_r)g
$$
that holds in $C_{\psi_2,\infty}(X)$ for any $g \in C_{\psi_1}$ according to (i).
 Passing to the adjoint operators it implies
$$
\frac{d}{dr} (g, V^{r,s} Y) = ((A_r - B_r)g, V^{r,s}Y)
$$
showing that $V^{r,s}$ yield a solution to \eqref{eqA14}. To show the uniqueness we shall
 use the method for the reduction of
the uniqueness problem to the existence of certain solutions of the adjoint problem,
 see e.g. \cite{Ma2} in the Hilbert space
setting and time independent generators. Let $0 < a < b < r$, $\chi_{[a,b]}(s)$ be
 an indicator function of $[a,b]$, and
 $v \in C_{\psi_1}(X)$. As $U^{t,r}$ solve \eqref{eqA13} one deduces that the function
$$
\phi_t = \int_t^r U^{s,r} \chi_{[a,b]}(s) v \, ds
$$
 solves the problem
\begin{equation}\label{eqA17}
\frac{d}{dt} \phi_t = - (A_t - B_t) \phi_t + \chi_{[a,b]}(s) v, \quad \phi_r = 0,
\end{equation}
in the sense that $\phi_t$ is continuous and satisfies \eqref{eqA17} everywhere
with possible exception of two points, where its derivative is not continuous.
Now, to prove uniqueness for\eqref{eqA14} it is enough to show that its any solution
 with $\xi_0 = 0$ vanishes. Assume that $\xi_t$ is a weakly continuous function
  in $\MC_{\psi_2}(X)$ such that $\xi_0 = 0$ and \eqref{eqA14} holds for all $g \in C_{\psi_1}$.
   Integration by parts, \eqref{eqA14} and weak continuity of $\xi_t$ imply that
$$
0 = (\phi_t, \xi_t)\mid^r_{t=0} = \int_0^r [(\dot{\phi}_t, \xi_t) + ((A_t - B_t) \phi_t, \xi_t)]dt
$$
whenever $\phi_t$ has a uniformly bounded derivatives in $C_{\psi_2,\infty}(X)$ apart from
 a finite number of points. Using
\eqref{eqA17} yield now the equation
$$
\int_b^a (v, \xi_t)\, dt = 0.
$$
 As it holds for arbitrary $0 < a < b < r$, $v \in C_{\psi_1}(X)$, it implies that $\xi_t = 0$.

(iii) From \eqref{eqA14} it follows that
\begin{equation}\label{eqA18}
(g, \xi_r) - (g, \xi_s) = \int_s^r ((A_t - B_t) g, \xi_t)\, dt, \quad 0 \leq s \leq r.
\end{equation}
Approximating any $g \in B_{\psi_1}(X)$ by functions from $C_{\psi_1}$ and using the dominated
 convergence one concludes that
\eqref{eqA18} holds for $g \in B_{\psi_1}(X)$. From this one deduces that $\xi_t$
 is an absolutely continuous function of $t$
in the norm $\MC_{\psi_1}(X)$. From boundedness of $\xi_t$ in $\MC_{\psi_2 f}(X)$ (that follows from weak continuity)
 and the weak continuity in $\MC_{\psi_1}(X)$ it follows the continuity in $\MC_{\psi_2 f}(X)$.
At last, again from \eqref{eqA18} one concludes that $\xi_t$ is continuously differentiable in $\MC_{\psi_1 f}(X)$.

(iv) This is  straightforward. Namely, one compares $U^{r,s}$ for various
 $\omega$ by a formula similar to \eqref{eqA12}. This yields the continuous
  dependence of $U^{r,s} \phi$ on $\omega$ for $\phi \in C^{\psi_1}(X)$.
   By approximation one extends this result to all $\phi \in C^{\psi_2}_\infty$.

\paragraph{Acknowledgments.} The author is grateful to S. Gaubert and M. Akian for their hospitality in INRIA (France),
 the major part of this work being done during the author's visit to INRIA in the spring 2006
  and to Martine Verneuille for the excellent
  typing of this manuscript in Latex.

\end{document}